\documentclass[onecolumn]{autart}
\usepackage{pdfpages}

\usepackage{algorithmic,url,color}
\usepackage{graphicx}
\usepackage{lscape}
\usepackage{algorithm}
\usepackage{amssymb,cite}
\usepackage{amsmath,amsfonts,mathrsfs}
\newtheorem{Lemma}{Lemma}
\newtheorem{Corollary}{Corollary}

\newtheorem{Definition}{Definition}
\newtheorem{Remark}{Remark}
\newtheorem{Notation}{Notation}
\newtheorem{Theorem}{Theorem}
\newtheorem{Example}{Example}

\newtheorem{Problem}{Problem}
\newcommand{\Rank}{\mathrm{Rank}}

\newcommand{\V}{\mathcal{V}}
\newcommand{\IM}{\mathrm{im}}
\newcommand{\LMAP}{\mathcal{M}}

\newenvironment{IEEEproof}{\begin{pf}}{ \rule{1ex}{1ex} \end{pf}}

\begin{document}
\begin{frontmatter}
\title{Solutions of differential-algebraic equations as outputs of LTI systems: application to LQ control problems}
\author{Mih\'aly Petreczky and } \ead{mihaly.petreczky@mines-douai.fr}
\author{Sergiy Zhuk}
\ead{sergiy.zhuk@ie.ibm.com}
\address{
%Ecole des Mines de Douai, Douai, France \\
CNRS, Centrale Lille, UMR 9189 - CRIStAL-Centre de Recherche en Informatique, Signal et Automatique de Lille, F-59000 Lille, France \\
IBM Research, Dublin, Ireland}

\maketitle
\begin{abstract}
 In this paper we synthesize behavioral ideas with geometric control theory and propose a unified geometric framework for representing all solutions of a Linear Time Invariant Differential-Algebraic Equation (DAE-LTI) as outputs of classical Linear Time Invariant systems (ODE-LTI). An algorithm for computing an ODE-LTI that generates solutions of a given DAE-LTI is described. It is shown that two different ODE-LTIs which represent the same DAE-LTI are feedback equivalent. The proposed framework is then used to solve an LQ optimal control problem for DAE-LTIs with rectangular matrices.
% and (ii) to realize the solutions of these problems as outputs of LTI systems.
\end{abstract}
\end{frontmatter}

\section{Introduction}
\label{sect:intro}
Consider a linear time invariant differential-algebraic equation (abbreviated by DAE-LTI) of the form
 \begin{equation}
   \label{eq:DAE}
     \dfrac {d(Ex(t))}{dt}=Ax(t)+Bu(t)
 \end{equation}
with arbitrary rectangular matrices $E,A \in \mathbb{R}^{c \times n}$ and $B \in \mathbb{R}^{c \times m}$. In this paper we discuss how to represent solutions of~\eqref{eq:DAE} as outputs of linear time invariant ordinary differential equations (abbreviated by ODE-LTI). This representation is then applied to derive necessary and sufficient solvability conditions for LQ optimal control problems with DAE-LTI constraints.

%The need for such a framework arises in several applications, where the DAE-LTIs at hand are not-regular and the solution is not smooth and cannot be interpreted as distributions.
Non-regular DAE-LTIs in the form~\eqref{eq:DAE} arise in control from several sources. They could either be a result of modeling physical systems, or arise
as a result of interconnecting several (possibly regular) DAEs. Indeed, regular DAE-LTIs are not closed under interconnection and so by applying a state-feedback to a regular DAE-LTIs one may arrive at a non-regular DAE-LTI~\cite{StefanovskiTAC}. Another nice example of non-regular DAE-LTIs are ODE-LTI with unknown external inputs. For instance, such systems arise when approximating Partial Differential Equations (PDEs) by ODEs. Then the approximation error can be viewed as an unknown input~\cite{ZhukCDC13,SZJFIHBS13IEEETAC}. Such systems can be modelled by DAE-LTIs if the inputs are viewed as a part of the state.
LQ optimal control for DAE-LTIs in the form~\eqref{eq:DAE} was studied by many authors~\cite{LQGNonReg,Stefanovski1,Stefanovski1.1,Stefanovski2,Kurina1,Kurina2,Kurina3,BergerThesis,DAEBookChapter}. In the mentioned papers solutions of~\eqref{eq:DAE} were defined either as smooth functions or as distributions. In many applications, however, the solution of~\eqref{eq:DAE} cannot be assumed to be smooth or to be a distribution, the former being too restrictive, while the latter does not correspond to the physical meaning of the DAE-LTI state. The present paper is motivated by the need for a framework which (A) provides a simple description of all solutions of DAE~\eqref{eq:DAE}, for which $Ex$ is absolutely continuous and $x,u$ are locally integrable, and
(B) allows to efficiently compute solutions of LQ optimal control problems for such DAE-LTIs. As an example of an application which requires such a framework, we mention the problem of state estimation for DAE-LTIs which arises in numerical analysis.\\
Many linear PDEs can be viewed as a linear time-invariant system $u_t = \mathcal{A}u$, $u(0)=u_0\in V$ with an infinite dimensional state space $V$ (for instance, Sobolev spaces of weakly differentiable functions). The precise choice of $V$ depends on the type of a differential operator $\mathcal{A}$. Often, it is possible to find a suitable system of orthogonal basis vectors in $V$ and identify $u$ with the infinite vector of its coordinates w.r.t. to this basis. In order to compute the solution $u$, the infinite vector $u=(u_1,u_2,\ldots,u_k,\ldots)^T$ is approximated by its truncation $u^h=\mathcal P^h u=(u_1,u_2,\ldots,u_h)^T$. In many applications, part of the state $u$ can be measured experimentally, i.e. finite dimensional measurement vectors $y(t)=\mathcal{C}u(t)$ are available, where $\mathcal{C}$ is an ``observation operator''. In ~\cite{ZhukCDC13,SZJFIHBS13IEEETAC} it was shown
that the truncated vector $z:=u^h$ satisfies a DAE-LTI in the following form:
\begin{equation}
\label{eq:dae_output}
   \begin{split}
      & \dfrac{d(Fz(t))}{dt} =Gz(t)+f\,, \\
      & y(t) = Hz(t) +\eta(t)\,,
  \end{split}
 \end{equation}
 where $F,G,H$ are certain rectangular matrices representing truncations of $\mathcal{A}$ and $\mathcal{C}$, and $f$, $\eta$ represent the terms which model the effect of the truncation error. Note that the time derivative of $z=u^h$ is a function of all the components of $u=(u_1,\ldots,u_k,\ldots)^T$, not just the first $h$ ones, hence the error terms $f$ and $\eta$. In \cite{ZhukCDC13,SZJFIHBS13IEEETAC} it was shown that for certain classes of PDEs and certain choice of basis functions,
%it ispossible to show that
%possible to estimate the energy ($L^2$ norm) of $f,\eta$, assuming that $u(0)$ belongs to
% a bounded subset of $V$. Intuitively, it means that the initial states belong to some ball around the origin, but the ball is now interpreted in a functional space. In this case, it can be shown
$(z(0),f,\eta)$ belong to the set
 %$x_0\in\mathbb R^n$, $f(t)\in\mathbb R^m$ and $\eta\in\mathbb R^p$ are uncertain and belong to an ellipsoid:
 \begin{equation*}
% \label{eq:dae_output:ellips}
   \begin{split}
    & \mathscr E=\{(z_0,f,\eta) \mid \rho(Fz_0,f,\eta)\le 1\}\,,\\
    &\rho(q_0,f,\eta):=q^T_0\hat{Q}_0q_0+
     \int_0^{\infty} (f^TQ^{-1}f+\eta^TR^{-1}\eta) dt\,.
   \end{split}
 \end{equation*}
 for suitable positive definite matrices $\hat{Q}_0,Q, R$. That is, \cite{ZhukCDC13,SZJFIHBS13IEEETAC} proposes a Galerkin-style method for solving PDEs, but unlike the classical methods, it takes into account the truncation error explicitly. The problem is that the obtained equation \eqref{eq:dae_output} cannot be solved numerically, since $f$ and $\eta$ are not known. However, one could use the experimental data $y(t)$ to estimate
 $z(t)$. Since we have bounds on the norms of $\eta$, $f$ and $z(0)$,  we could use a minimax observer to find an estimate of $z(t)$ such that the maximal (worst-case) difference between this estimate and $z(t)$
 is minimal. %Here, minimization is performed on the set of potential observers.
From~\cite{ZhukPetrezky2013TAC,ZhukPetreczkyCDC13} it follows that in order to construct such an observer, we have to solve an LQ control problem, that is to minimize $\rho$ over solutions of a dual DAE-LTI given by~\eqref{eq:DAE} with $A=-G^T,B=H^T,E=F^T$.

 Notice that according to~\cite{ZhukCDC13,SZJFIHBS13IEEETAC} the state $z$ of~\eqref{eq:dae_output} is absolutely continuous as it models $u^h$, and $F$ is a rectangular matrix. Hence, (i) the dual DAE-LTI will not be regular and may have several solutions (or none at all) from any initial state and any input (see Example \ref{interconnect}), and (ii) each solution of the dual DAE-LTI will have an absolutely continuous part and a measurable part. Thus, the usual assumption on regularity, impulse controllability, etc. do not hold for the dual system. Moreover, the solutions cannot be
 assumed either to be smooth or to be a distribution. These observations clearly indicate that a framework satisfying conditions (A) and (B) is required to estimate the state of~\eqref{eq:dae_output} and obtain a robust approximation of PDE's solution.

The aim of the present paper is to propose a framework featuring (A) and (B) for general DAEs-LTI. To this end we use behavioral approach~\cite{WillemsBehavior} and geometric control theory~\cite{TrentelmanBook}: namely, given DAE-LTI in the form~\eqref{eq:DAE} we achieve point (A) above by introducing a class of associated ODE-LTIs:
\begin{equation}
  \label{eq:DAELin}
  \begin{split}
    &\dot p = A_lp+B_lv  \\
    & w  =  C_lp+D_lv
  \end{split}
\end{equation}
such that the external behavior (set of output trajectories) of \eqref{eq:DAELin} coincides with the set of solutions of the given DAE-LTI, and~\eqref{eq:DAELin} satisfies a number of nice technical conditions detailed in the following section. %A generic procedure which constructs the entire class of associated ODE-LTIs represents the first contribution of this paper.

Representing solutions of DAE-LTIs as outputs of ODE-LTIs is a classical idea. The earliest method relies on Kronecker canonical form \cite{Gantmacher,DAEBookChapter}. However, the method requires to differentiate the inputs and so either the input is assumed to be smooth or the solution is viewed as a distribution. We take another approach which is based on the observation that solutions of a DAE-LTI can be viewed
as output nulling solutions of a suitable ODE-LTI. Various versions of this approach appeared in \cite{Zhuk2012sysid,LQGNonReg,Stefanovski1,Stefanovski2,Zhang2013}. However, in the cited papers, ODE-LTIs played a role of an auxiliary tool, and hence the constructions were not general, but rather problem specific and existence conditions were tailored to meet the requirements of the problem at hand. In contrast, this paper
describes the entire class of associated ODE-LTIs which have a simple system-theoretic interpretation: they are
feedback equivalent to a  minimal ODE-LTI realizations in the sense of \cite{WillemsBehavior} of the solution set of the DAE-LTI at hand.
The results \cite{Zhuk2012sysid,LQGNonReg,Zhang2013} are special cases of the ones which are presented here, provided the corresponding assumptions are used. The construction of
\cite{Stefanovski1,Stefanovski2} is closely related, but it is not formally a special case due to the different solution concept used in the paper.
%There is a long history of geometric ideas in DAEs, see for example \cite{Malabre,Armentano} and the references in the surveys \cite{LewisGeomTut,DAEBookChapter}. However, to the best of our knowledge, these ideas were not used to represent all absolutely continuous solutions of DAE~\eqref{eq:DAE} as outputs of ODE-LTIs. We stress that relating DAEs with ODEs via Weierstrass canonical forms is a classical trick. However, this trick requires either smoothness constraints on admissible inputs or usage of distributions as it was mentioned above. We note that despite the superficial difference, there is a connection between our approach and that of based on quasi-Weistrass canonical forms \cite{TrennWeierstrass}. We postpone the detailed discussion until Remark \ref{weierstrass:rem1}.

The concept of ODE-LTIs associated with DAE-LTIs allows us to easily achieve point (B) above, namely solve the infinite horizon LQ control problem for DAE-LTIs, by reducing it to the classical LQ control problem for ODE-LTIs.
%%given positive definite matrices $Q,R$ and positive semi-definite matrix $Q_0$, and for a given initial state $x_0$, find a solution $(x^{*},u^{*})$ of \eqref{eq:DAE} with $Ex^*(0)=Ex_0$, such that \(J_{\infty}(x^{*},u^{*})=\limsup_{t \rightarrow \infty} \inf_{(x,u)} J_t(x,u) \), where the infimum is taken over all the solutions $(x,u)$ of \eqref{eq:DAE} such that $Ex(0)=Ex_0$, and the cost functions $J_\infty$, $J_t$ are defined as follows:
%%   \begin{equation}
%%    \label{eq:DAE:lqr}
%%      \begin{split}
%%       & J_{\infty}(x,u)=\limsup_{t \rightarrow \infty}  J_t(x,u)\,. \\
%%       & J_t(x,u) =(x^T(t)TE^TQ_0Ex(t) +  \\
%%      & + \int_0^{t} x^T(s)Qx(s)ds + u^T(s)Ru(s)ds\,.
%%      \end{split}
%%   \end{equation}
In particular, we derive new necessary and sufficient solvability conditions for the infinite horizon LQ control problem in terms of behavioral stabilizability of DAE-LTIs. Specifically, we show that
% we show that the LQ control problem has a solution from
%some initial state $x_0$ if and only if $Ex_0$ belongs to the behavioral stabilizability subspace of~\eqref{eq:DAE} (see definition~\ref{def:stabil} in Section~\ref{lqr:control}).
the optimal value of the quadratic cost function is given by a norm of the initial condition which is induced by a unique solution of the algebraic Riccati equation. Moreover, if $(x^*,u^*)$ is the optimal trajectory verifying $Ex^*(0)=Ex_0$, then
$u^{*}=Kx^{*}$ for some matrix $K$, i.e. the optimal input has the form of a feedback.
Note that this does not imply that all solutions of the closed-loop system $\dfrac{d(Ex)}{dt}=(A+BK)x$, $Ex(0)=Ex_0$ are optimal, as the latter may have several solutions, including ones which render the cost function infinite.
% In other words, the closed loop DAE~\eqref{eq:DAE} may have solutions $x$ for which the cost is not minimal or even infinite.
However, we show that there exist matrices $K_1,K_2$ such that $(x^*,u^*)$ is the only solution of the DAE-LTI
$\dfrac{d(Ex)}{dt}=Ax+Bu, K_1x+K_2u=0$ which satisfies $Ex(0)=Ex_0$. The additional algebraic constraint $K_1x+K_2u=0$ can be thought of as a generalization of state-feedback concept and can be interpreted as a
controller in the sense of behavioral approach~\cite{WillemsBehavior} (see Remark \ref{interconnect}).
Note that controllers which are not of feedback form can still be implemented and in fact are widely used for controlling physical devices \cite{WillemsInterconnect}.
Moreover, for the purposes of observer design~\cite{ZhukPetrezky2013TAC,ZhukPetreczkyCDC13} it is sufficient that at least one trajectory of the closed-loop system is optimal. %% as there, it is sufficient to find a stable autonomous ODE-LTI whose output is $(x^*,u^*)$.

The literature on optimal control for DAE-LTIs is vast. For an overview we refer to~\cite{LewisGeomTut} and the references therein. To the best of our knowledge, the most relevant references on LQ control are \cite{BenderLaub1987IEEETAC,LQGFrank,LQGNonReg,Stefanovski1,Stefanovski1.1, Stefanovski2,Kurina1,Kurina2,Kurina3}.
In \cite{BenderLaub1987IEEETAC,LQGFrank} only regular DAEs were considered.
The infinite horizon LQ control problem for non-regular DAE was also addressed in \cite{LQGNonReg}, however there it is assumed that the DAE has a solution from any initial state. %In terms of results, \cite{Stefanovski1,Stefanovski1.1,SWith respect to \cite{Stefanovsk1,Stefanovski2} the main differences are as follows.
We consider existence of a solution from a particular initial condition, as opposed to \cite{Stefanovski1,Stefanovski1.1, Stefanovski2}.
This is done both for the sake of generality and in order to address the requirements of already mentioned observer design problems~\cite{ZhukPetreczkyCDC13,ZhukPetrezky2013TAC}.
 %: note that if we want to estimate the state component $\ell^TFz(t)$, then we are interested only in those  optimal solutions of the dual system which start at $F^T\ell^T$.
Furthermore, in contrast to \cite{Stefanovski1,Stefanovski1.1, Stefanovski2,Stefanovski4}, where only sufficient conditions
are presented, in this paper we present necessary and sufficient conditions.
Moreover, the cost function considered here differs from the one of \cite{Stefanovski1,Stefanovski1.1,Stefanovski2,Stefanovski4}, as it includes
a terminal cost term $x^T(t_1)E^TQ_0Ex(t_1)$. Note that the latter term is indispensable to transform an observer design problem into a dual control problem (see~\cite[Theorem 1]{ZhukPetrezky2013TAC} for the further details). In addition, we allow non-smooth solutions. This leads to subtle but important technical differences.\\
We note that in \cite{BruellLQDAEBehav} the behavioral approach was used for LQ control of DAE-LTIs, however, the LQ problem considered in
\cite{BruellLQDAEBehav} is different from the one of this paper and it does not present detailed algorithms.
%Furthermore,  we impose different additional restrictions on the optimal
%solution. As a result, we can weaken some of the conditions imposed in~\cite{Stefanovski1,Stefanovski1.1,Stefanovski2,Stefanovski4}.
%\item
%\item
% Finally, similarly to \cite{Stefanovski1,Stefanovski2}, we solve the LQ problem
%by reducing it to an LTI LQ problem. Unlike in
%\cite{Stefanovski1,Stefanovski2}, where Silverman's algorithm was
%used, we obtain corresponding LTI through geometric arguments. The
%construction of \cite{Stefanovski1,Stefanovski2} then becomes a
%particular case of our construction, if we restrict $x$ to be only absolutely continuous.
%\end{enumeradrate}
%Optimal control of non-linear and time-varying DAEs was also addressed in the literature, see \cite{Kurina2007,MehramnnKunk2008,KunkelMehrmann1997}. However, they do not seem to be directly applicable to the problem of the current paper, as they look at finite horizon optimal control problems.
The results of \cite{Kurina1,Kurina2,Kurina3,Kurina4} provide sufficient conditions for existence of an optimal controller for stationary DAEs: these conditions involve existence of a solution to an algebraic Riccati equation. In contrast, we provide conditions which are necessary and sufficient, and are, therefore, less restrictive. To illustrate this we describe an LQ control problem for a simple DAE-LTI such that the
conditions of \cite{Kurina1,Kurina2,Kurina3,Kurina4} are not satisfied (see discussion after Example \ref{interconnect}).
This LQ control problem arises as a dual of an observer design problem for a DAE-LTI of the form~\eqref{eq:dae_output}. On the other hand, this generality comes at price: the sufficient conditions of \cite{Kurina1,Kurina2,Kurina3,Kurina4} yield a feedback such that all trajectories of the closed-loop system are optimal. In contrast, the solution of this paper does not always yield such a feedback law.

An extended version of this paper is available at
\cite{AutomaticaPaperArxive} and its preliminary version appeared in \cite{ZhukPetreczkyCDC13}. With respect to
\cite{ZhukPetreczkyCDC13} the main difference is that we included detailed proofs, and provided necessary and sufficient conditions for
existence of a solution for the infinite horizon optimal control
problem. The solution of the finite horizon optimal control problem
was already presented in~\cite{Zhuk2012sysid}. In contrast
to~\cite{Zhuk2012sysid} we consider the infinite horizon case
too, and the algorithm of~\cite{Zhuk2012sysid} for computing an ODE-LTI that generates solutions of a given DAE-LTI is one of many possible implementations of the generic procedure of this paper.

\textbf{Outline of the paper}
In Section \ref{dae:geo} we present the notion of an ODE-LTI associated with a DAE-LTI and prove that all ODE-LTI representing the same DAE-LTI are feed-back equivalent. In Section \ref{lqr:control} we apply this result to solve the infinite horizon LQ control problem.
In Section \ref{sect:num} we present a numerical example. 

\textbf{Notation}
$I_n$ denotes the $n \times n$ identity matrix; for an $n \times n$ matrix $S$,
$S>0$ means $x^TSx>0$ for all $0 \ne x\in\mathbb R^n$,  %, $\|S\|$ denotes the Frobenius norm;
$F^+$ denotes the Moore-Penrose pseudoinverse of the matrix $F$.
Consider an interval $I \subseteq \mathbb{R}$ of the form $[a,b]$, $[a,+\infty)$, $(-\infty,a]$, $a,b \in \mathbb{R}$, or $I=\mathbb{R}$.
For an integer $p > 0$ denote by $L^p(I,\mathbb{R}^{n})$ (or simply by $L^p(I)$) the space of all measurable functions $f:I \rightarrow \mathbb{R}^{n}$ such that $\int_I ||f||^{p} dm < +\infty$, where $m$ is the Lebesgue measure on $\mathbb{R}$ (see \cite{Rudin:CompReal} for more details). Let $L^p_{loc}(I,\mathbb{R}^{n})=\{f:I \rightarrow \mathbb{R}^n \mid f \in L^p(K,\mathbb{R}^{n}),\forall K \subsetneq I, K \mbox{ is a compact interval } \}$,
i.e. the restriction of $f\in L^p_{loc}(I,\mathbb{R}^n)$ onto any compact sub-interval $K$ of $I$ is in $L^p(K,\mathbb{R}^n)$. Note that $L^p_{loc}(I,\mathbb{R}^n)=L^p(I,\mathbb{R}^n)$ for compact intervals  $I$. We will use the usual conventions to denote integrals with respect to the Lebesgue measure, see \cite[page 52, Remark 2.21]{Rudin:CompReal}. In particular, $\int_I f dm$, $\int_a^b f(s)ds$ denote the same integral for $I=[a,b]$, $a,b \in \mathbb{R}$ or $I=[a,+\infty)$.
%$int_a^b f(s)ds$, if $I=[a,+\infty)$, $a \in \mathbb{R}$, the we use $\int_a^{+\infty} f(s)ds$. Sometimes, in order to avoid  notation overload, we will use other symbols than $s$ and $ds$ in the expression
%inside the integral, but in all these case we still mean Lebesgue integral on a sub-interval of $\mathbb{R}$.
Denote by $AC(I,\mathbb{R}^{n})$ the set of all absolutely continuous functions $f:I \rightarrow \mathbb{R}^{n}$, see \cite{Rudin:CompReal} for the definition of absolute continuity. Note that if $f \in AC(I,\mathbb{R}^n)$, then there exists a function $g\in L^p_{loc}(I,\mathbb{R}^n)$ such that $f(t)=f(0)+\int_0^{t} g(s)ds$, $\forall t \in I$.
In accordance with the convention, \cite{Rudin:CompReal}, we say that an equation $f_1(t)=f_2(t)$ holds almost everywhere (write $f_1(t)=f_2(t)$ a.e. or simply $f_1=f_2$ a.e.) for any measurable functions $f_1,f_2:I \rightarrow \mathbb{R}$, if there exists a set $S \subseteq I$, such that $S$ is of Lebesgue measure zero and
for any $t \in I, t \notin S$, $f_1(t)=f_2(t)$. Finally, $f|_{A}$ stands for the restriction of a function $f$ onto a set $A$.

\section{Linear systems associated with DAEs}
\label{dae:geo}
Consider a linear time-invariant differential-algebraic system (DAE-LTI)
 \begin{equation}
\label{dae:sys}
\dfrac{dEx(t)}{dt} = Ax(t) + Bu(t)\,.
\end{equation}
Here $A,E \in \mathbb{R}^{c \times n}$, $B \in \mathbb{R}^{c \times m}$. In this section we will define a class of ODE-LTIs whose output trajectories are the state and input trajectories of \eqref{dae:sys} and show that these ODE-LTIs exist and they are unique up to feedback equivalence. To this end, we view the set of solutions of \eqref{dae:sys} as behaviors in the sense of \cite{WillemsBehavior,WillemsLinear}, and we view the ODE-LTIs as their state-space representations. We then state a number of consequences of this fact for the solvability theory of \eqref{dae:sys}. The section is organized as follows. In \S \ref{dae:geo:main} we present the main results. In \S \ref{dae:geo:proof} we present the proofs of the results.

\subsection{Main results}
\label{dae:geo:main}
In order to carry out the program outlined above, we start by defining solutions of \eqref{dae:sys}.
\emph{In this section, by an interval we mean an interval of one of the following forms:
$I=\mathbb{R}$, $I=[a,+\infty)$, $I=(-\infty,a]$, $I=[a,b]$, $a,b \in \mathbb{R}$.}
\begin{Definition}
\label{sol:not}
\label{sol:behav}
 Let $I$ be an interval. %$\mathbb{R}$, or $[0,+\infty)$ or $[0,t]$, $0 < t \in \mathbb{R}$.
 A \emph{solution} of \eqref{dae:sys}  on the interval $I$ is a tuple $(x,u) \in L^1_{loc}(I,\mathbb{R}^n) \times L^1_{loc}(I,\mathbb{R}^m)$ such that
$Ex$ is absolutely continuous, and $\forall t_0 \in I$, $\forall t \in I, t \ge t_0: Ex(t)=Ex(t_0)+\int_{t_0}^{t} (Ax(s)+Bu(s))ds$.
%$(x,u)$ satisfy \eqref{dae:sys} almost everywhere.
%$\dfrac{dEx}{dt} = Ax(t) + Bu(t)$ a.e..
%%\end{Definition}
%  Let $t_0 \in I$ be an arbitrary element of the connected interval $I$.
% Note that $(x,u) \in L^1_{loc}(I,\mathbb{R}^n) \times L^1_{loc}(I,\mathbb{R}^m)$ is a solution of \eqref{dae:sys} on $I$ if and only if $Ex$ is absolutely continuous and
% for all $t \in I$, $Ex(t)=Ex(t_0)+\int_{t_0}^{t} (Ax(s)+Bu(s))ds$. The latter means that solutions can be defined using integral equations.
 %That is, the condition that the differential equation \eqref{dae:sys} holds almost all $t \in I$ is equivalent to the condition that
 %the integral equation $Ex(t)=Ex(t_0)+\int_{t_0}^{t} (Ax(s)+Bu(s))ds$ holds for all $t \in I$.
%\begin{Definition}[Behavior $\mathcal{B}_I(E,A,B)$]
Denote by $\mathcal{B}_{I}(E,A,B)$ the set of all solutions $(x,u) \in L^1_{loc}(I,\mathbb{R}^n) \times L^1_{loc}(I,\mathbb{R}^m)$ of \eqref{dae:sys} defined
on $I$.
%The set $\mathcal{B}_{I}(E,A,B)$ will be referred to as the
%\emph{behavior} of \eqref{dae:sys} on $I$ or simply behavior.
\end{Definition}
That is, $(x,u) \in L^1_{loc}(I,\mathbb{R}^n)  \times L^1_{loc}(I,\mathbb{R}^m)$ is a solution of \eqref{dae:sys} on $I$ if and only if $Ex$ is absolutely continuous and
$\dfrac{dEx(t)}{dt} = Ax(t) + Bu(t)$ a.e.
%%In the sequel, we will also need the set of all pairs of $(x,u)  \in L^1_{loc}(I,\mathbb{R}^n) \times L^1_{loc}(I,\mathbb{R}^m)$, whose modification on a set of measure zero belongs to $\mathcal{B}_I(E,A,B)$.
%\begin{Definition}[Extended behavior $\mathcal{B}_{I,e}(E,A,B)$]
% With the notation of Definition \ref{sol:behav},
% let $\mathcal{B}_{I,e}(E,A,B)=\{ (x,u) \in L^1_{loc}(I,\mathbb{R}^n) \times L^1_{loc}(I,\mathbb{R}^m) \mid \exists (\hat{x},\hat{u}) \in \mathcal{B}_I(E,A,B), \hat{x}=x, \hat{u}=u \mbox{ a.e. } \}$.
%\end{Definition}
%That is, if $(x,u) \in \mathcal{B}_{I,e}(E,A,B)$, then $(x,u)$ is not necessarily a solution, but it can become one after we modify $x$ and $u$ on a set of measure zero. In particular,
%$\mathcal{B}_I(E,A,B) \subseteq \mathcal{B}_{I,e}(E,A,B)$.
%Note that in our definition of a solution we follows what could be called a behavioral approach.
Note that solutions may happen to be non-smooth or even discontinuous (except $Ex$), so they may contain jumps. Also distributions (as DAE's solutions) are not allowed. Hence, in our setting the solution of DAE-LTI has no ``impulsive parts''. We stress that if we allowed for distributional solutions then DAE-LTI would have solutions with impulsive parts as we do not restrict matrices $E,A,B$.

Let us now recall few definitions from the behavioral approach \cite{WillemsLinear}. Consider a linear time-invariant system defined by differential equation (referred as ODE-LTI),
 \begin{equation}
 \label{ode:lti}
     \begin{split}
      & \dot p = \widetilde{A}p+Gq\,,\\
      & z = \widetilde{C}p+\widetilde{D}q,
     \end{split}
 \end{equation}
 where $\widetilde{A} \in \mathbb{R}^{r \times r}$, $G \in \mathbb{R}^{r \times s}$, $\widetilde{C} \in \mathbb{R}^{p \times r}$,  $\widetilde{D} \in \mathbb{R}^{p \times s}$. We identify the ODE-LTI \eqref{ode:lti} with the corresponding tuples $(\widetilde{A},G,\widetilde{C},\widetilde{D})$ of matrices. Let $I$ be an interval and $\mathcal{B} \subseteq  L^1_{loc}(I,\mathbb{R}^p)$. Following the definition of \cite{WillemsLinear}, we say that the ODE-LTI \eqref{ode:lti} is a \emph{realization} of $\mathcal{B}$, if
 \textbf{(1)} for every $z \in \mathcal{B}$ there exist functions $p \in AC(I,\mathbb{R}^{r})$,
 $q \in L^1_{loc}(I,\mathbb{R}^s)$ such that $\dot p=\widetilde{A}p+Gq$ a.e., and $z=\widetilde{C}p+\widetilde{D}q$ a.e., and
 \textbf{(2)} if $(p,q) \in AC(I,\mathbb{R}^{r}) \times  L^1_{loc}(I,\mathbb{R}^s)$  is such that  $\dot p=\widetilde{A}p+Gq$ a.e. then $z=\widetilde{C}p+\widetilde{D}q$ a.e. for some $z \in \mathcal{B}$.
  %then there exists $\widetilde{z} \in \mathcal{B}$ such that
  %$z=\widetilde{z}$ a.e.
  That is, if \eqref{ode:lti} is a realization of $\mathcal{B}$, then any element of $\mathcal{B}$ is an output
  trajectory of \eqref{ode:lti},  and conversely any output trajectory of \eqref{ode:lti} belongs to $\mathcal{B}$, possibly after having been modified on a set of measure zero. In the sequel, we are interested in ODE-LTI realizations of $\mathcal{B}_I(E,A,B)$. Note that $\mathcal{B}_I(E,A,B)$ can naturally be viewed as a subset of $L^1_{loc}(I,\mathbb{R}^{n+m})$, so the
 definition above can be applied.  With this terminology we define the notion of a ODE-LTI  system associated with a DAE-LTI.
\begin{Definition}
\label{def:asscoLTI}
\label{linassc:def}
  An ODE-LTI system of the form
   \begin{equation}
       \label{asscoLTI}
        \mathscr{S}\left\{\begin{split}
         & \dot v = A_lv+B_lg\,,  \\
         & \nu = C_{l}v+D_{l}g\,. \\
       \end{split}\right.
  \end{equation}
  $A_l \in \mathbb{R}^{\hat{n} \times \hat{n}}$, $B_l \in \mathbb{R}^{\hat{n} \times k}$, $C_l \in \mathbb{R}^{(n+m) \times \hat{n}}$,
  $D_l \in \mathbb{R}^{(n+m) \times k}$, $\hat{n} \le n$,
  is called an ODE-LTI associated with the DAE-LTI~\eqref{dae:sys}, if the following conditions hold:
  \begin{enumerate}
\item Either $\begin{bmatrix} D^T_l, & B^T_l \end{bmatrix}^T=0$ and $k=1$, or $D_l$ is full column rank: $\Rank D_l=k$.
\item
  Let $C_s$ and $D_s$ be the matrices formed by the first
  $n$ rows of $C_l$ and $D_l$ respectively. Then $ED_s=0$, $\Rank EC_s = \hat{n}$.
 \item
    $\mathscr{S}=(A_l,B_l,C_l,D_l)$ is a realization of $\mathcal{B}_{I}(E,A,B)$
    for any interval $I$.
 \end{enumerate}
\end{Definition}
\begin{Notation}[$\LMAP$]
With the notation above, $\LMAP$ denotes the Moore-Penrose inverse of $EC_s$. The matrix $\LMAP$ will be referred to as the state map of
 $(A_l,B_l,C_l,D_l)$.
\end{Notation}
\begin{Theorem}
\label{dae2lin:theo_main}
Let $(A_l,B_l,C_l,D_l)$ be an ODE-LTI associated with \eqref{dae:sys}. If
$(x,u) \in \mathcal{B}_I(E,A,B)$, and we define  the function $v = \LMAP Ex$ and $g=D_l^{+}((x^T,u^T)^T-C_l\LMAP Ex)$, then
$v \in AC(I,\mathbb{R}^{\hat{n}})$, $g \in L_{loc}^1(I,\mathbb{R}^k)$, and
\begin{equation}
\label{dae2lin:theo_main:eq1}
     \dot v = A_lv+B_lg \mbox{ a.e. } \quad (x^T,u^T)^T = C_lv + D_lg \mbox{ a.e. .}
\end{equation}
%Moreover, for any $(v,g) \in AC(I,\mathbb{R}^{\hat{n}}) \times L_{loc}^1(I,\mathbb{R}^k)$ which satisfies \eqref{dae2lin:theo_main:eq1},
 Conversely,  for any $(v,g) \in AC(I,\mathbb{R}^{\hat{n}}) \times L_{loc}^1(I,\mathbb{R}^k)$ such that
 $\dot v = A_lv+B_lg \mbox{ a.e }$, $C_lv+D_lg \in \mathcal{B}_I(E,A,B)$.
\end{Theorem}
That is, not only the outputs of the associated ODE-LTI correspond to the solutions of the DAE-LTI, but the state trajectory of the DAE-LTI determines the corresponding state trajectory of the associated ODE-LTI.
%To conclude, below we show that an associated linear system exists for every DAE.

The question arises if associated ODE-LTIs exist. The answer is affirmative.
\begin{Theorem}[Existence]
\label{dae2lin:theo}
 Consider the DAE-LTI system \eqref{dae:sys}.
 There exists an  ODE-LTI system $\mathscr{S}=(A_l,B_l,C_l,D_l)$ associated with \eqref{dae:sys}.
\end{Theorem}
The proof of Theorem \ref{dae2lin:theo} is constructive and it yields an easy to implement algorithm for computing an associated ODE-LTI.
The Matlab code of the algorithm is available at \verb'http://sites.google.com/site/mihalypetreczky/'.

The next question is whether associated ODE-LTIs of the same DAE-LTI are related in any way. In order to answer this question we need the following terminology.
\begin{Definition}[Feedback equivalence]
  Two ODE-LTIs $\mathscr{S}_i=(A_i,B_i,C_i,D_i)$, $i=1,2$ are said to be
  \emph{feedback equivalent}, if there exist
  a matrix $K$ and two nonsingular square matrices $U,T$ of suitable dimensions such that
  $(T(A_1+B_1K)T^{-1},TB_1U,(C_1+D_1K)T^{-1},D_1U)=\mathscr{S}_2$.
 We will call $(T,K,U)$ feedback equivalence
 from $\mathscr{S}_1$ to $\mathscr{S}_2$.
\end{Definition}
\begin{Theorem}
\label{dae2lin:theo:uniq}
 Any two ODE-LTI systems associated with the same DAE-LTI \eqref{dae:sys} are feedback equivalent.
\end{Theorem}

 Existence and uniqueness of associated ODE-LTIs allow us to study existence of solutions for DAE-LTIs from a given initial state.
% It is know that general DAE-LTIs do not admit a solution for
% any boundary condition $Ex(0)=z$. This prompts us to introduce the following terminology.
 \begin{Definition}[Consistency set $\mathcal{V}(E,A,B)$]
\label{sol:consist}
  We will say that a vector $z \in \mathbb{R}^{c}$ is \emph{differentiably consistent}, if
  there exists a solution $(x,u)$ of \eqref{dae:sys} defined on an interval $I \subseteq \mathbb{R}$ such that $0 \in I$ and $Ex(0)=z$. We denote by
  $\mathcal{V}(E,A,B)$ the set of all differentiably consistent vectors $z \in \mathbb{R}^c$.
\end{Definition}
 %From linearity of \eqref{dae:sys} it follows that $\mathcal{V}(E,A,B)$ is a linear subspace of $\mathbb{R}^c$.
% Note that our definition of the set $\mathcal{V}(E,A,B)$ is similar to the definition of the set $\mathcal{V}^{diff}_{[E,A,B]}$ of all consistent initial differentiable variable of \cite[Section 2]{DAEBookChapter}, in fact, the two definitions can be shown to be equivalent.
 %%  except that the solution concept considered in
 %%\cite{DAEBookChapter} is slightly different. However, later on we will show that $\mathcal{V}^{diff}_{[E,A,B]}= \mathcal{V}(E,A,B)$.
  % Theorem \ref{dae2lin:theo} then yields the following corollaries.
   \begin{Corollary}
 \label{col1:important}
  %$z \in \mathbb{R}^c$ is is differentiably consistent if and only if $z \in \IM EC_s$. In particular,
  Let $(A_l,B_l,C_l,D_l)$ be an ODE-LTI associated with the DAE-LTI \eqref{dae:sys}, and let $C_s$ be the
  matrix formed by the first $n$ rows of $C_l$.  Then $\mathcal{V}(E,A,B)=\IM EC_s$.
 \end{Corollary}
 \begin{Corollary}
 \label{col2:important}
   Let $I$ be any interval such that $0 \in I$.
   If $z$ is differentiably  consistent,  there exists a solution $(x,u)$ of \eqref{dae:sys} on  $I$,
   such that $Ex(0)=z$. Moreover, $(x,u)$ can be chosen so that $x,u$ are smooth functions.
  \end{Corollary}
 In principle, it could happen that there exists a solution $(x,u)$ on the interval $[0,t_1]$ such that $z=Ex(0)$, but there exist no solution $(x,u)$ with $z=Ex(0)$ for a larger interval $[0,t_2]$, $t_2 > t_1$. In this case, the subsequent formulation of the finite and infinite horizon control problem would be more involved.
Corollary \ref{col2:important} tells us that this can never happen.
 Corollary \ref{col2:important} also implies that if
 there exist a solution $(x,u)$ on $I$ such that $Ex(0)=Ex_0$, then
 there exists a solution $(x,u)$ on $I$ with $Ex(0)=Ex_0$ and $x$ being differentiable, for any interval $I$ containing $0$.
   Finally, recall from \cite{DAEBookChapter} the notion of impulse controllability. Using
  \cite[Corollary 4.3]{DAEBookChapter}, we can show the following.
 \begin{Corollary}
 \label{col3:important}
 For any $x_0 \in \mathbb{R}^n$, $Ex_0$ is differentiably consistent $\iff$ \eqref{dae:sys} is impulse controllable $\iff$ for any matrix $Z$ such that $\IM Z = \ker E$,
 $\mathrm{rank} \begin{bmatrix} E, & A, & B \end{bmatrix} = \mathrm{rank} \begin{bmatrix} E, & AZ, & B \end{bmatrix}$.
 \end{Corollary}

 To conclude this section, we would like to discuss the relationship between the results above and existing results.
 To begin with, existence of an associated ODE-LTI is not that surprising.
 Note that \eqref{dae:sys} can be viewed as a kernel representation of $\mathcal{B}_{I}(E,A,B)$,
 if one disregards the subtleties related to smoothness of solutions. It is a classical result
 that behaviors admitting a kernel representation can be represented as outputs of ODE-LTIs \cite{WillemsBehavior,WillemsLinear}. What makes a separate proof of Theorem \ref{dae2lin:theo} necessary are the subtle issues related to
 differentiability of solutions and the additional properties we require for associated ODE-LTIs.
 In fact, the proof of Theorem \ref{dae2lin:theo} bears a close resemblance to \cite{Schumacher19881} which provides an algorithm for computing a state-space realization of a kernel representation of a behavior.
 However, unlike \cite{Schumacher19881}, the proof of Theorem \ref{dae2lin:theo} exploits the specific
 structure of DAE-LTIs and yields existence of ODE-LTI realizations which satisfy Definition  \ref{def:asscoLTI}.
% applies only to smooth behaviors and the obtained ODE-LTI does not necessarily satisfy the conditions of

 Feedback equivalence of associated ODE-LTIs stems from minimality theory for behaviors.
 Recall that according to \cite{WillemsLinear} an ODE-LTI $(A,B,C,D)$ is a minimal realization of a  behavior
 $\mathcal{B} \subseteq L^1_{loc}(\mathbb{R},\mathbb{R}^p)$, if $(A,B,C,D)$ is a realization of $\mathcal{B}$ and
 for any other ODE-LTI $(A^{'},B^{'},C^{'},D^{'})$ which is a realization of $\mathcal{B}$,
  the number of state variables of
 $(A,B,C,D)$ is not greater than the number of state variables of $(A^{'},B^{'},C^{'},D^{'})$. % and the number of inputs of $(A,B,C,D)$ is not greater than the number of inputs of $(A^{'},B^{'},C^{'},D^{'})$.
% We can tacitly restrict attention to minimal realizations $(A,B,C,D)$
% otherwise we can replace the input space of $(A,B,C,D)$ while keeping it a realization of $\mathcal{B}$.
 In \cite{WillemsLinear} it was shown that any two minimal state-space representations of the same
 behavior are feedback equivalent. It turns out that associated ODE-LTIs are in fact minimal:

 \begin{Corollary}[Minimality]
\label{col:minim}
  %Assume that there exist at least one ODE-LTI system  associated with \eqref{dae:sys}. Then the following holds.
 \begin{enumerate}
 \item Any ODE-LTI system associated with \eqref{dae:sys} is a minimal state-space representation of the behavior
       $\mathcal{B}_{\mathbb{R}}(E,A,B)$.
 \item Conversely, if $(A,B,C,D)$  is a minimal realization of  $\mathcal{B}_{\mathbb{R}}(E,A,B)$
 such that $\begin{bmatrix} B^T & D^T \end{bmatrix}^T$ is  either full column rank or it is zero and has one column,
 then $(A,B,C,D)$ is an ODE-LTI associated with \eqref{dae:sys}.
 %\item Any two ODE-LTI systems associated with the same DAE-LTI \eqref{dae:sys} are feedback equivalent.
 \end{enumerate}
\end{Corollary}
\subsection{Proofs}
\label{dae:geo:proof}
\begin{IEEEproof}[Theorem \ref{dae2lin:theo_main}]
 Let $C_s,D_s$ be the matrices formed by the first $n$ rows of $C_l$, $D_l$. From Definition \ref{def:asscoLTI} it follows that $ED_s=0$ and $EC_s$ is full column rank.
 Let $(x,u) \in  \mathcal{B}_{I}(E,A,B)$. Since $(A_l,B_l,C_l,D_l)$ is a realization of $\mathcal{B}_I(E,A,B)$, it follows that there exist
 $(v,g) \in AC(I,\mathbb{R}^n) \times L_{loc}^1(I,\mathbb{R}^k)$ such that \eqref{dae2lin:theo_main:eq1} holds. It then follows that
 $Ex=EC_sv+ED_sg=EC_sv$ a.e., since $ED_s=0$. Note that $Ex$ and $EC_sv$ are both absolutely continuous, hence $Ex=EC_sv$ a.e. implies that
 $Ex(t)=EC_sv(t)$ \emph{for all} $t \in I$.
 Finally, $(x^T,u^T)=C_lv+D_lg$ a.e. and the fact that $D_l$ is either zero or it is full column rank,
 imply that $g=D_l^{+}((x^T,u^T)^T-C_lv)=D_l^{+}((x^T,u^T)^T-C_l\LMAP Ex)$ a.e.
 In order to show the second statement, notice that since  $(A_l,B_l,C_l,D_l)$ is a realization of $\mathcal{B}_{I}(E,A,B)$, there exist $(x,u) \in \mathcal{B}_I(E,A,B)$ such that
 $(x^T,u^T)^T = C_lv+D_lg$ a.e. Let $(\widetilde{x}^T,\widetilde{u}^T) = C_lv+D_lg$ with $\widetilde{x} \in L^1_{loc}(I,\mathbb{R}^n)$ and $\widetilde{u} \in  L^1_{loc}(I,\mathbb{R}^m)$.
 It then follows that $E\widetilde{x}=EC_sv+ED_sg=EC_sv$ is absolutely continuous. Since $x=\widetilde{x}$ a.e. and $u=\widetilde{u}$ a.e., absolute continuity of $E\widetilde{x}$ implies
 that $(\widetilde{x},\widetilde{u}) \in \mathcal{B}_I(E,A,B)$.
\end{IEEEproof}

In order to present the proof of Theorem \ref{dae2lin:theo}, we recall the following notions from geometric control theory of linear systems.  Consider an ODE-LTI
of the form \eqref{ode:lti}.
 Let $I$ be an interval and let $t_0 \in I$.
 %Let $u \in L^1_{loc}(I,\mathbb{R}^{q})$ and let $x_0 \in \mathbb{R}^{r}$.
 %Let $x(x_0,t_0,u) \in AC(I,\mathbb{R}^r)$ and $y(x_0,t_0,u) \in L^1_{loc}(I,\mathbb{R}^p)$  be the state and output trajectory of
 %\eqref{ode:lti} which corresponds to the initial state $x_0=x(t_0)$ and input $u$, i.e. $x(x_0,t_0,u)(t_0)=x_0$
 %$\dot x(x_0,u)= Ax(x_0,t_0,u)+Bu$ a.e. and $y(x_0,t_0,u)=Cx(x_0,t_0,u)+Du$ a.e.
 Recall from \cite[Definition 7.8]{TrentelmanBook}
 the concept
 of a weakly unobservable subspace of the ODE-LTI \eqref{ode:lti}.
 I.e.,
  an initial state $p_0 \in \mathbb{R}^{r}$ of \eqref{ode:lti}
  is \emph{weakly unobservable}, if there exist $p \in AC([0,+\infty),\mathbb{R}^r)$ and  $q \in L^1_{loc}([0,+\infty), \mathbb{R}^s)$ such that
  $\dot p=\widetilde{A}p+Gq$ a.e., $p(0)=p_0$,  $0=\widetilde{C}p+\widetilde{D}q$ a.e.
   Following \cite{TrentelmanBook},
 let us denote the set of all
 weakly unobservable states by $\V$. Recall from \cite[Section 7.3]{TrentelmanBook}, $\V$ is a vector space and in fact
 it can be computed.
 For technical purposes we will need the following easy extension of \cite[Theorems 7.10--7.11]{TrentelmanBook}.
\begin{Theorem} %{\cite[Theorem 7.10--.11]{TrentelmanBook}}
\label{theo:geo_contr}
 Consider the ODE-LTI \eqref{ode:lti}. With the notation above:
\begin{enumerate}
%\item \highlight{Empty space????}
 %\label{theo:geo_contr:part1}
 %     $\V=\V(\LIN)$ is the largest subspace of $\mathbb{R}^{r}$
 %     for which it holds that
 %     $\begin{bmatrix} \widetilde{A} \\ \widetilde{C} \end{bmatrix} \V \subseteq \V \times 0 + \IM \begin{bmatrix} G \\ \widetilde{D} \end{bmatrix}$.
\item
 \label{theo:geo_contr:part2}
$\V$ is the largest subspace of $\mathbb{R}^{r}$ for which
      there exists a linear map
      $\widetilde{F}:\mathbb{R}^r \rightarrow \mathbb{R}^{s}$ such that
      \begin{equation}
      \label{theo:geo_contr:eq1}
        (\widetilde{A}+G\widetilde{F})\V \subseteq \V \mbox{ and }
        (\widetilde{C}+\widetilde{D}\widetilde{F})\V = 0
      \end{equation}
\item
 \label{theo:geo_contr:part3}
     Let $\widetilde{F}$ be a map such that \eqref{theo:geo_contr:eq1} holds for
      $\V$. Let
      $L \in \mathbb{R}^{q \times k}$ for some $k$
      be a matrix such that $\IM L = \ker \widetilde{D} \cap G^{-1}(\V)$. Choose $L$ so that
      $L$ is full column rank if $\ker \widetilde{D} \cap G^{-1}(\V) \ne \{0\}$, or $L=0 \in \mathbb{R}^{q \times 1}$ otherwise.

      For any interval $I$, for any $p \in AC(I,\mathbb{R}^r)$, $q \in L^1_{loc}(I,\mathbb{R}^s)$,
      \[ \widetilde{C}p+\widetilde{D}q=0  \mbox{ for } t \in I \mbox{ a.e.} \]
       if and only if $p(t) \in \V$ for all $t \in I$, and there exists $g \in L^1_{loc}(I,\mathbb{R}^{k})$
       such that:
      \[ q(t)=\widetilde{F}p(t)+Lg(t) \mbox{ for } t \in I \mbox{ a.e.} \]
\end{enumerate}
\end{Theorem}
\begin{IEEEproof}[Proof of Theorem \ref{theo:geo_contr}]
 Part \ref{theo:geo_contr:part2} is a
 reformulation of \cite[Theorem 7.10]{TrentelmanBook}.
 For $I=[0,+\infty)$, Part \ref{theo:geo_contr:part3} is a restatement of
 \cite[Theorem 7.11]{TrentelmanBook}. For all the other intervals $I$,  the proof is similar to
 \cite[Theorem 7.11]{TrentelmanBook}.
\end{IEEEproof}
\begin{IEEEproof}[Proof of Theorem \ref{dae2lin:theo}]
 There exist suitable nonsingular matrices $S$ and $T$ such that
 \begin{equation}
 \label{dae:sys:tr}
     SET = \begin{bmatrix} I_r & 0 \\
                          0   & 0
          \end{bmatrix},
 \end{equation}
where $r = \Rank ~E$. Let
 \begin{equation}
 \label{dae:sys:tr1}
   SAT=\begin{bmatrix} \widetilde{A} & A_{12} \\
                    A_{21} & A_{22}
     \end{bmatrix} \mbox{,\ \ }
     SB=\begin{bmatrix} B_1 \\ B_2 \end{bmatrix}
   %  HT=\begin{bmatrix} H^T_1 \\ H_2^T \end{bmatrix}^T
 \end{equation}
 be the decomposition of $E,A,B$ such that $\widetilde{A} \in \mathbb{R}^{r \times r}$ and $B_{1} \in \mathbb{R}^{r \times m}$.
 Define
  \begin{equation}
 \label{dae:sys:tr2}
   \begin{split}
    & G = \begin{bmatrix} A_{12}, & B_1 \end{bmatrix}  \mbox{, \ \ }
     \widetilde{D} = \begin{bmatrix} A_{22}, & B_2 \end{bmatrix}  \mbox{ and }
     \widetilde{C}=A_{21}
   \end{split}.
 \end{equation}
 %Then \eqref{dae:sys1} is equivalent to
 %\begin{equation}
 %\label{dae:sys2}
  %\begin{split}
  %   & \dot p = \widetilde{A}p+Gq  \widetilde{C}p + \widetilde{D}q
  %  \end{split}
 %\end{equation}
  Consider the ODE-LTI \eqref{ode:lti} with the choice of $\widetilde{A},G,\widetilde{C},\widetilde{D}$ as defined in \eqref{dae:sys:tr1} and \eqref{dae:sys:tr2}.
  %with the state $p$, input $q$ and output $z$
% \begin{equation}
% \label{dae:sys_lin}
%    \begin{split}
%     & \dot p = \widetilde{A}p+Gq  \\
%     &  z     = \widetilde{C}p + \widetilde{D}q
%    \end{split}.
% \end{equation}
  We claim that for any $(x,u) \in L^1_{loc}(I,\mathbb{R}^n) \times L^1_{loc}(I,\mathbb{R}^{m})$,  $(x,u) \in \mathcal{B}_{I}(E,A,B)$ if and only if,
  the functions $(p,q) \in L^1_{loc}(I,\mathbb{R}^r) \times L^1_{loc}(I,\mathbb{R}^{n+m-r})$ defined by
  \begin{equation}
  \label{dae2lin:theo:pf:eq1}
    \begin{bmatrix} p \\ q \end{bmatrix} =\begin{bmatrix} T^{-1} & 0 \\ 0 & I_m \end{bmatrix}\begin{bmatrix} x \\ u \end{bmatrix}
  \end{equation}
  are such that $p \in AC(I,\mathbb{R}^r)$, $\dot p = \widetilde{A}p+Gq$ a.e. and $0=\widetilde{C}p+\widetilde{D}q$ holds for all a.e..
  Indeed, notice that  $SEx=SETT^{-1}x=(p^T,0)^T$.
  Hence, $Ex$ is absolutely continuous if and only if $p$ is absolutely continuous. Furthermore, notice that
  \[ S(Ax+Bu)=\begin{bmatrix} \widetilde{A}p+Gq \\ \widetilde{C}p+\widetilde{D}q \end{bmatrix}. \]
  Hence, $\dfrac{dEx}{dt} = Ax+Bu$ a.e. if and only if
  \[ \begin{bmatrix} \dot p \\ 0 \end{bmatrix} = S \dfrac{dEx}{dt} = S(Ax+Bu)= \begin{bmatrix} \widetilde{A}p+Gq \\ \widetilde{C}p+\widetilde{D}q \end{bmatrix} \mbox{ a.e. } \]
  Finally note that $(x,u) \in \mathcal{B}_I(E,A,B)$ if and only if $Ex$ is absolutely continuous and $\dfrac{dEx}{dt} = Ax+Bu$ a.e..
The desired linear system $\mathscr{S}=(A_l,B_l,C_l,D_l)$ may be obtained as follows.
 Let $\widetilde{F}$ and $L$ be the matrices from Theorem \ref{theo:geo_contr} applied to the ODE-LTI $(\widetilde{A},G,\widetilde{C},\widetilde{D})$, and let $\V$ be the space of
 weakly unobservable initial states of $(\widetilde{A},G,\widetilde{C},\widetilde{D})$.
  %In particular, $(\widetilde{A}+\widetilde{F}G)\V \subseteq \V$, $\IM L=\ker \widetilde{D} \cap G^{-1}(\V)$
 %and $L$ is either full column rank (if $\ker \widetilde{D} \cap G^{-1}(\V) \ne \{0\}$ or $L=0 \in \mathbb{R}^{(n-r+m) \times 1}$ otherwise.
 Define the matrices
 \[
 \begin{split}
 & \bar{C} = \begin{bmatrix} T & 0 \\ 0 & I_m \end{bmatrix} \begin{bmatrix} I_r \\ \widetilde{F} \end{bmatrix} \mbox{ and }
  \bar{D} = \begin{bmatrix} T & 0 \\ 0 & I_m \end{bmatrix}\begin{bmatrix} 0 \\ L \end{bmatrix}.
 \end{split}
\]
 From Theorem \ref{theo:geo_contr}  and the discussion above it then follows that for any
 $(p,g) \in AC(I,\mathbb{R}^r) \times L^1_{loc}(I,\mathbb{R}^{k})$  such that $p \in AC(I,\mathbb{R}^r)$, $p(t) \in \V$, $t \in I$,  $\dot p =(\widetilde{A}+G\widetilde{F})p+GLg$ a.e., it holds that
 $\widetilde{C}p+\widetilde{D}(\widetilde{F}p+Lg))=0$ and hence
 $(x^T,u^T) = \bar{C}p+\bar{D}g$ belongs to $\mathcal{B}_I(E,A,B)$. Conversely,
 if $(x,u) \in \mathcal{B}_{I}(E,A,B)$, then  there exist
 $(p,g) \in AC(I,\mathbb{R}^r) \times L^1_{loc}(I,\mathbb{R}^{k})$  such that $p \in AC(I,\mathbb{R}^r)$, $p(t) \in \V$, $t \in I$,  $\dot p =(\widetilde{A}+G\widetilde{F})p+GLg$ a.e. and
 $(x^T,u^T)^T = \bar{C}p+\bar{D}g$ a.e. .
%Consider the decomposition $\bar{C}_{s}:\V \rightarrow \mathbb{R}^n$,
%$\bar{C}_{inp}:\V \rightarrow \mathbb{R}^m$,
%$\bar{D}_{s}:\mathbb{R}^{k} \rightarrow \mathbb{R}^n$
%$\bar{D}_{inp}:\mathbb{R}^{k} \rightarrow \mathbb{R}^m$
%such that $\bar{C}(p)=\begin{bmatrix} \bar{C}_s(p)^T, & \bar{C}^T_{inp}(p)\end{bmatrix}^T$
%and $\bar{D}(p)=\begin{bmatrix} \bar{D}_s(v)^T & \bar{D}^T_{inp}(v)\end{bmatrix}^T$.
%It is then easy to see that $Ex_0 \in \mathcal{X}$ if and only
%if $T^{-1}x_0=(p^T,0)^T$ for some $p \in \V$.
 %Moreover,  it can be verified that $\IM E\bar{C}_s=\mathcal{X}$ and
%that $\begin{bmatrix} I_r & 0 \end{bmatrix} S$ is the left inverse of $E\bar{C}_s$ and thus $\Rank E\bar{C}_s=\hat{n}$.
 %Furthermore, it can easily be checked that $E\bar{D}_s=0$.
 %$\bar{C}^T\bar{D}=0$ and $\bar{D}^T\bar{D}=I_k$.
 Consider a basis $b_1,\ldots,b_r$ of $\mathbb{R}^r$ such that
 $b_1,\ldots,b_{\hat{n}}$ span $\V$. Let $\mathscr{R}  \in \mathbb{R}^{r \times r}$ be the  corresponding basis transformation, i.e. $\mathscr{R}^{-1}=\begin{bmatrix} b_1,& \ldots, & b_r \end{bmatrix}$.
 Let $D_l=\bar{D}$ and $A_l,B_l,C_l$ be the matrix representations of the linear maps $(\widetilde{A}+G\widetilde{F})|_{\V}:\V \rightarrow \V$, $\bar{C}|_{\V}:\V \rightarrow \mathbb{R}^{n + m}$, $GL:\mathbb{R}^k \rightarrow \V$
 respectively in the basis $b_1,\ldots, b_{\hat{n}}$. That is,
 %$\mathscr{R}(\V)=\IM \begin{bmatrix} I_{\hat{n}} & 0 \\ 0 & 0 \end{bmatrix}$.
 \[
   \begin{split}
   & A_l=\begin{bmatrix} I_{\hat{n}} & 0  \end{bmatrix} \mathscr{R} (\widetilde{A}+G\widetilde{F}) \mathscr{R}^{-1} \begin{bmatrix} I_{\hat{n}} \\ 0  \end{bmatrix}, \quad
  C_l =  \bar{C}  \mathscr{R}^{-1} \begin{bmatrix} I_{\hat{n}} \\ 0  \end{bmatrix}, \\
  & B_l= \begin{bmatrix} I_{\hat{n}} & 0  \end{bmatrix} \mathscr{R} GL.
  \end{split}
 \]
%% Choose a basis of $\V=\V(\LIN)$ and choose $(A_l,B_l,C_l,D_l)$ as follows:
%% $D_l=\bar{D}$, and let $A_l,B_l,C_l$ be the matrix representations in this basis of
%% the linear maps $(\widetilde{A}+G\widetilde{F}):\V \rightarrow \V$, $GL:\mathbb{R}^k \rightarrow \V$, and
%% $\bar{C}:\V \rightarrow \mathbb{R}^{n+m}$ respectively.
% Define
% \[ \mathcal{X}=\{ ET\begin{bmatrix} p \\ 0 \end{bmatrix} \mid p \in \V\}. \]
 It is easy to see that  with this choice,  $\mathscr{S}=(A_l,B_l,C_l,D_l)$ satisfies
 Definition \ref{def:asscoLTI}.
\end{IEEEproof}
%\begin{Remark}[Uniqueness of the linear system associated with DAE]
%\label{rem:uniq}
%% The proof of Theorem \ref{dae2lin:theo} is constructive and yields
%% an algorithm for computing $(A_l,B_l,C_l,D_l)$ from $(E,A,B)$.
%% This prompts us to introduce the following terminology.
%%\begin{Definition}
%%\label{linassc:def}
%% A linear system $\mathscr{S}=(A_l,B_l,C_l,D_l)$ described in the proof of
%% Theorem \ref{dae2lin:theo} is called
%% the linear system associated with the DAE \eqref{dae:sys}, and the map
%% $\LMAP$ is called the corresponding state map.
%%\end{Definition}
\begin{Remark}
\label{geo:rem1}
  Notice that the dimension $\hat{n}$ of the associated linear system constructed in the proof of Theorem \ref{dae2lin:theo}
  satisfies $\hat{n} \le \Rank E \le \max\{c,n\}$.
\end{Remark}
\begin{Remark}[Comparison with~\cite{Zhuk2012sysid}]
 The system $(A_1^s,A_2^s)$ described in~\cite[Proposition 3]{Zhuk2012sysid} is related
 to the ODE-LTI constructed in the proof of Theorem \ref{dae2lin:theo}, see \cite{AutomaticaPaperArxive} for a detailed explanation.

\end{Remark}

\begin{Remark}%[Relationship to Wang sequence and quasi-Weierstrass canonical forms]
 \label{weierstrass:rem1}
  Recall from \cite{DAEBookChapter} that the augmented Wong sequence is
  defined as follows $\V_0=\mathbb{R}^n$, $\V_{i+1}=A^{-1}(E\V_i+\IM B)$, and
  that the limit $\V^*=\bigcap_{i=0}^{\infty} \V_i$ is achieved in a finite
  number of steps: $\V^*=\V_k$ for some $k \in \mathbb{N}$.
  It is not difficult to see that $\V$ from the proof of
  Theorem \ref{dae2lin:theo} correspond to the limit $\V_{*}$ of the
  augmented Wong sequence $\V_i$ for the DAE \eqref{dae:sys}:
  $\V^{*}=\{ (p,Fp+Lq)^T \mid p \in \V\}$.
  Hence, if $(A_l,B_l,C_l,D_l)$ is an ODE-LTI associated with \eqref{dae:sys}, then
  $\V^{*}=\IM \begin{bmatrix} C_s & D_s \end{bmatrix}$, where $C_s,D_s$ are the matrices formed by the first
  $n$ rows of $C_l$ and $D_l$ respectively.
  In\cite{TrennWeierstrass} a relationship between
  the quasi-Weierstrass form of regular DAEs and space $\V^{*}$ for $B=0$ was established. This indicates that there
  might be a deeper connection between quasi-Weierstrass forms and
  associated linear systems. The precise relationship remains a topic of future research.
 \end{Remark}
 Before presenting the proof of Theorem \ref{dae2lin:theo:uniq}, we present the proof of
 Corollary \ref{col:minim}, as it yields Theorem \ref{dae2lin:theo:uniq}.
\begin{IEEEproof}[Proof of Corollary \ref{col:minim}]
   Consider an ODE-LTI system $(A_l,B_l,C_l,D_l)$ associated with \eqref{dae:sys}. From \cite[Theorem 4.3]{WillemsLinear} it follows $(A_l,B_l,C_l,D_l)$ is a minimal realization of
  $\mathcal{B}_{\mathbb{R}}(E,A,B)$, if and only if $\mathcal{V}^{*}=0$, where
  $\mathcal{V}^{*}$ is the set of weakly unobservable states of $(A_l,B_l,C_l,D_l)$.
  %=\{ p_0 \in \mathbb{R}^{\hat{n}} \mid \exists (q,p) \in L_{loc}^1(\mathbb{R},\mathbb{R}^k) \times AC(\mathbb{R},\mathbb{R}^{\hat{n}}): p(0)=p_0, \quad \dot p(t) = A_lp(t)+B_lq(t), \mbox{ a.e. } \quad C_lp(t)+D_lq(t)=0 \mbox{ a.e. }\}$.
 %For any $(q,p) \in L_{loc}^1(\mathbb{R},\mathbb{R}^k) \times AC(\mathbb{R},\mathbb{R}^{\hat{n}})$,$C_lp(t)+D_lq(t)=0$ a.e.
Let $C_s,D_s$ be the matrices formed by the first $n$ rows of $C_l$ and respectively $D_l$.
For any $p_0 \in \mathcal{V}^{*}$, there exist $(q,p) \in L_{loc}^1(\mathbb{R},\mathbb{R}^k) \times AC(\mathbb{R},\mathbb{R}^{\hat{n}}): p(0)=p_0$, $\dot p(t) = A_lp(t)+B_lq(t), \quad C_lp(t)+D_lq(t)=0$ a.e..
This implies $E(C_sp(t)+D_sq(t))=EC_sp(t)=0$ a.e.  and by continuity of $p(t)$ this implies $EC_sp(t)=0$, $t \in I$.
In particular, $p_0=p(0)=(EC_s)^{+}(C_sp(0))=0$. That is, $\mathcal{V}^{*}=0$ and $(A_l,B_l,C_l,D_l)$ is minimal.

Let $(A,B,C,D)$ be a minimal realization of $\mathcal{B}_{\mathbb{R}}(E,A,B)$, such that either $\begin{bmatrix} B \\ D \end{bmatrix}$ is full column rank or it is zero.
Let $(A_l,B_l,C_l,D_l)$ be an ODE-LTI system associated with \eqref{dae:sys}. From the discussion above it follows that $(A_l,B_l,C_l,D_l)$ is a minimal realization of  $\mathcal{B}_{\mathbb{R}}(E,A,B)$.
From the proof of \cite[Theorem 7.1]{WillemsLinear} it follows that there exists a nonsingular matrix $T$ and a matrix $K$ such that
$A_l=T(A_l+B_lK)T^{-1}$, $\IM \begin{bmatrix} TB_l, \\ D_l \end{bmatrix} = \IM \begin{bmatrix} B \\ D \end{bmatrix}$ and $(C_l+D_lK)T^{-1}=C$.
Hence, if $D_l=0,B_l=0$, then $B=0,D=0$ and both $\begin{bmatrix} TB_l, \\ D_l \end{bmatrix}$ and  $\begin{bmatrix} B \\ D \end{bmatrix}$ have one column. Thus, with $U=1$, $TB_lU=B$, $D_lU=D$ holds.
Otherwise, since both $\begin{bmatrix} TB_l \\ D_l \end{bmatrix}$ and  $\begin{bmatrix} B \\ D \end{bmatrix}$ are full column rank,  $\begin{bmatrix} TB_l \\ D_l \end{bmatrix}$ and  $\begin{bmatrix} B \\ D \end{bmatrix}$
have the same number of columns and there exists an invertible matrix $U$ such that $U=1$, $TB_lU=B$, $D_lU=D$.
Hence, $(A,B,C,D)$ is feedback equivalent with $(A_l,B_l,C_l,D_l)$. Using feedback equivalence, it is easy to see that $(A,B,C,D)$ satisfies all the conditions of Definition \ref{def:asscoLTI}.
\end{IEEEproof}
\begin{IEEEproof}[Theorem \ref{dae2lin:theo:uniq}]
 %From Corollary \ref{col:minim}, ODE-LTI systems associated with the same DAE-LTI \eqref{dae:sys} are minimal
 %realization of $\mathcal{B}_{\mathbb{R}}(E,A,B)$. From the proof of Corollary are feedback equivalent.
 Theorem \ref{dae2lin:theo:uniq} is a direct consequence of the proof of Corollary \ref{col:minim}.
\end{IEEEproof}
 \begin{IEEEproof}[Proof of Corollary \ref{col1:important}]
  If $z$ is differentiably consistent, then there exists a solution $(x,u)$ of \eqref{dae:sys} on $I$ for some interval $I$, such that $0 \in I$ and $Ex(0)=z$.
  Then there exist $(v,g) \in AC(I,\mathbb{R}^{\hat{n}}) \times L^1_{loc}(I,\mathbb{R}^k)$ such that $\dot v = A_lv+B_lg$ and $(x^T,u^T)^T=C_lv+D_lg$ a.e.. In particular, $Ex(t)=EC_sv(t)$ for almost all $t \in I$, and hence, by continuity of $Ex$ and $v$, $Ex(t)=EC_sv(t)$ for all $t \in I$, where $C_s$ is the matrix formed by the first $n$ rows of $C_l$.
  %$t_n \in I$ such that $|t_n| < 2^{-n}$ and $Ex(t_n)=EC_sv(t_n)$.
  % Note that $v$ and $Ex$ are both absolutely continuous and hence
  %they are continuous.
  Therefore, $z=Ex(0)=EC_sv(0)$ and hence
  $z \in \IM EC_s$. Conversely, if $z=EC_sv_0$, then let $v$ be the solution of $\dot v = A_lv$, $v(0)=v_0$  on $I$ and set
  $(x^T,u^T)^T=C_sv$. Then $Ex(0)=EC_sv(0)=EC_sv_0=z$ and $(x,u)$ is a solution of \eqref{dae:sys} on $I$, i.e. $z$ is differentiably consistent.
 \end{IEEEproof}
 % From te discussion above, we can deduce the following corollary.
   \begin{IEEEproof}[Proof of Corollary \ref{col2:important}]
    Consider an ODE-LTI $(A_l,B_l,C_l,D_l)$ which is an associated ODE-LTI for the DAE-LTI \eqref{dae:sys}.
    If $z$ is differentiably consistent, then $z \in \IM EC_s$. Let $p \in AC(I,\mathbb{R}^{\hat{n}})$ be the solution of the differential equation $\dot p = A_lp$, $p(0)=\LMAP(z)$. Then $p$ is smooth. From the properties of the associated ODE-LTI it then follows that $(x^T,u^T)^T=C_lp$ is a solution of \eqref{dae:sys} which satisfies $Ex(0)=z$.
    Moreover, as $x$ and $u$ are linear functions of $p$, they are also smooth.
  \end{IEEEproof}
  \begin{IEEEproof}[Proof of Corollary \ref{col3:important}]
Recall from~\cite[Section 2]{DAEBookChapter} that $\mathcal{V}^{diff}_{[E,A,B]}$ is the set of differentiably consistent initial conditions: $x_0 \in \mathcal{V}_{[E,A,B]}^{diff}$ if and only if there exists a solution $(x,u)$ of \eqref{dae:sys} on $\mathbb{R}$ such that $Ex(0)=Ex_0$ and $x$ is absolutely continuous. From Corollary \ref{col2:important} it follows that $\mathcal{V}(E,A,B)=\mathcal{V}^{diff}_{[E,A,B]}$. The statement follows now from \cite[Corollary 4.3]{DAEBookChapter}.
 \end{IEEEproof}

%Since it is known from \cite{WillemsLinear} that any behavior which admits a realization by an LTI admits a minimal LTI realization, Corollary \ref{col:minim} implies the following.
%That is, ODE-LTI systems associated with a DAE, if they exist, coincide with minimal linear state-space realizations of
%the behavior of the DAE, and thus they are unique up to feedback equivalence.

\section{Application to finite and infinite horizon LQ problem for DAEs}
\label{lqr:control}

 In this section we present the application of the results of Section \ref{dae:geo} to LQ control of DAE-LTIs. We will start by stating the problem formally.
 Consider a DAE-LTI of the form \eqref{dae:sys} We define the set of solutions which satisfy the boundary condition $Ex(0)=z$.
 \begin{Definition}%[Solutions from initial state $\mathscr{D}_{z}(t_1)$]
 Let $z \in \mathcal{V}(E,A,B)$.
 For $t_1 > 0$, denote by \( \mathscr{D}_{z}(t_1) \) the set of all solutions
\( (x,u) \)  of \eqref{dae:sys} on $[0,t_1]$ such that $Ex(0)=z$. % and $(x,u) \in L_2(I,\mathbb{R}^{n}) \times L_2(I,\mathbb{R}^m)$.
  Likewise, define   \( \mathscr{D}_z(\infty) \) as the set of all solutions $(x,u)$ of \eqref{dae:sys} on $[0,+\infty)$ such that $Ex(0)=z$.
  %and $(x,u) \in L^{loc}_2([0,+\infty),\mathbb{R}^{n}) \times L^{loc}_2([0,+\infty), \mathbb{R}^m)$.
  %as the set of all tuples $(x,u,d) \in A_E^{loc}(\mathbb{R}^n) \times L^2_{loc}(\mathbb{R}^{m}) \times L^2_{loc}(\mathbb{R}^n)$ such that
 %$(x,u)$ satisfies \eqref{dae:sys} and $Ed(t)=0$ for all $t \ge 0$.
 %A state $x_0$ will be called \emph{differentiably consistent} if $\mathscr{D}_{x_0}(t_1) \ne \emptyset$.
 %The set of differentiably consistent initial states is denoted by $\mathcal{V}_c$.
\end{Definition}
 %While we defined the set $\mathscr{D}_z(t_1)$ for any $z \in \mathbb{R}^{c}$, it is clear that this set can be non-empty only for $z \in \mathbb{R}^c$
 %of the form $z=Ex_0$, $x_0 \in \mathbb{R}^{n}$.
 From Corollary \ref{col2:important} it follows that $\mathscr{D}_{z}(t_1)$, $t_1 > 0$, $\mathscr{D}_{z}(\infty)$
 are not empty for any $z \in \mathcal{V}(E,A,B)$.

Take symmetric $R\in\mathbb R^{m\times m},Q,\in\mathbb R^{n\times n}$, $Q_0 \in \mathbb{R}^{c \times c}$ and assume that $R > 0, Q> 0$ $Q_0 \ge 0$.
 Fix an initial state $x_0 \in \mathbb{R}^n$.
 For any trajectory $(x,u) \in \mathscr{D}_{z}(t)$, $t \ge t_1$ define the cost
 functional
 \begin{equation}
 \label{opt1.1}
 \begin{split}
    & J_{t_1}(x,u) = x(t_1)^TE^TQ_0Ex(t_1)+\\
     & + \int_0^{t_1} (x^T(s)Qx(s)+u^T(s)Ru(s)) ds.
 \end{split}
 \end{equation}
  Note that in \eqref{opt1.1} $x,u$ may be defined on an interval larger than $[0,t_1]$, but $J_{t_1}(x,u)$ depends only on
  the restriction of $x$ and $u$ to $[0,t_1]$. Moreover, $J_{t_1}(x,u)$ need not be finite, as $x|_{[0,t_1]}$ and $u|_{[0,t_1]}$ may not belong to $L^2([0,t_1],\mathbb{R}^n)$ respectively $L^2([0,t_1],\mathbb{R}^m)$.
\begin{Problem}[Finite-horizon optimal control]
\label{opt:fin}
 Consider a differentiably consistent initial state $z \in \mathcal{V}(E,A,B)$.
 The problem of finding $(x^*,u^*) \in \mathscr{D}_{z}(t_1)$
 such that:
 \[ J_{t_1}(x^*,u^*) = J_{t_1}^{*} \stackrel{def}{=} \inf_{(x,u) \in \mathscr{D}_{z}(t_1)} J_{t_1}(x,u) < +\infty \]
 is called the \emph{finite-horizon optimal control problem for the initial state $z$} and
 $(x^*,u^*)$ is called the solution of the finite-horizon optimal control problem.
\end{Problem}
 Clearly, the optimal solution $(x^*,u^*)$ should be square integrable  (i.e. should belong to  $L^2([0,t_1],\mathbb{R}^n) \times L^2([0,t_1],\mathbb{R}^m)$)
 and when calculating $J^*_{t_1}$, infimum should be taken only over  $L^2([0,t_1],\mathbb{R}^n) \times L^2([0,t_1],\mathbb{R}^m) \cap  \mathscr{D}_{z}(t_1)$,
 since for all other solutions the cost function is infinite.
 \begin{Problem}[Infinite horizon optimal control]
\label{opt:contr:def}
 Consider a differentiably consistent initial state $z \in \mathcal{V}(E,A,B)$.
  For every $(x,u) \in \mathscr{D}_{z}(\infty)$, define
 \[ J_{\infty} (x,u)=\limsup_{t \rightarrow \infty} J_{t} (x,u)\,. \]
  The \emph{infinite horizon optimal control problem for the initial
  state $z$}  is the problem of finding $(x^*,u^*)$ such that $(x^{*},u^{*}) \in \mathscr{D}_{z}(\infty)$ and
  \begin{equation}
  \label{opt:eq2}
  \begin{split}
     & J_{\infty}(x^{*},u^{*}) =  J_{\infty}^{*}  < +\infty \\
     & J_{\infty}^{*}  \stackrel{def}{=} \limsup_{t_1 \rightarrow \infty } \inf_{(x,u) \in \mathscr{D}_{z}(t_1)} J_{t_1}(x,u)\,.
  \end{split}
 \end{equation}
  %and $(x^{*},u^{*})$ is the output of an autonomous stable linear system, i.e.
  %there exists matrices
 The pair $(x^{*},u^{*})$ will be called the solution of
 the infinite horizon (optimal) control problem for the initial state $z$.
\end{Problem}
 Note that the optimal solution $(x^*,u^*)$ of the infinite horizon problem should belong to  $L^2([0,+\infty),\mathbb{R}^n) \times L^2([0,+\infty),\mathbb{R}^m)$.
 Note that $L^2([0,+\infty),\mathbb{R}^p)$ is a subset of  $L^1_{loc}([0,+\infty),\mathbb{R}^p)$, $p=n,m$ and  hence
 $(x^*,u^*) \in L^2([0,+\infty),\mathbb{R}^n) \times L^2([0,+\infty),\mathbb{R}^m)$ does not conflict with the definition of
 $\mathscr{D}_{z}(\infty)$.
\begin{Remark}
 The proposed formulation of the infinite horizon control problem is not the most natural one. It also makes sense to look for solutions
 $(\tilde{x},\tilde{u}) \in \mathscr{D}_{z}(\infty)$ which satisfy
 \( J_{\infty} (\tilde{x},\tilde{u}) = \inf_{(x,u)\in \mathscr{D}_{z}(\infty)} J_{\infty} (x,u)\).
 The latter means that the cost induced by $(\tilde{x},\tilde{u})$ is the
 smallest among all the trajectories $(x,u)$ which are defined on the whole
 time axis.
 It is easy to see that if $(x^*,u^*)$ is a solution of  Problem \ref{opt:contr:def}, then
 \( J_{\infty} (x^{*},u^{*})=J_{\infty}^{*} = \inf_{(x,u) \in \mathscr{D}_{z}(\infty)} J_{\infty} (x,u) \), i.e. the solution of
 Problem \ref{opt:contr:def} yields the minimal cost among all the solutions $(x,u)$ of the DAE-LTI which satisfy $Ex=z$ and which are defined on the whole time axis.
 Another option is to use $\lim$ instead of $\limsup$ in the definition of $J_\infty(x^{*},u^{*})$ and in
 \eqref{opt:eq2}. In fact, the solution we are going to present remains valid if we replace $\limsup$ by $\lim$.
\end{Remark}
\begin{Remark}[Derivatives of inputs in $J$]
 Note that the cost function $J_{t_1}(x,u)$ does not contain explicitly the derivatives of $x$ and $u$. However, as it is well known from
 solution theory of DAEs, derivatives of $u$ can implicitly appear in the state $x$ and hence in the cost function. It is especially obvious if one computes the Kronecker canonical form of the DAE at hand and rewrites the cost function in the new coordinates. We stress that in our framework the state and input of the DAE-LTI are linear functions of the output of the associated ODE- LTI. Thus, if the DAE-LTI's state depends on derivatives of the input this will be taken into account implicitly. As a result, the cost $J$ will also include this relation implicitly (it is made clear in~\eqref{lpg:eq:cost1} where $J$ is reformulated in terms of the associated DAE-LTI).
\end{Remark}
The rest of the section is organized as follows. In \S \ref{lqr:control:main} we present the main results and in \S \ref{lqr:control:proof} we present their
proofs.

\subsection{Main results}
\label{lqr:control:main}

We start by presenting a solution to the
finite horizon case.
To this end, let $\mathscr{S}=(A_l,B_l,C_l,D_l)$ be an ODE-LTI associated with the DAE-LTI \eqref{dae:sys}.
  Consider the following differential Riccati equation
  \begin{equation}
  \label{DARE}
  \begin{split}
     & \dot P(t)= A_l^TP(t)+P(t)A_l-K^T(t)(D_l^TSD_l)K(t)+C_l^TSC_l \\
     & P(0)=(EC_s)^TQ_0EC_s, \quad S=\mathrm{diag}(Q,R) \\
     & K(t)=\left\{ \begin{array}{rl}
                (D_l^TSD_l)^{-1}(B_l^TP(t)+D_l^TSC_l) & \mbox{ if }  D_l\ne 0 \\
                0  \in \mathbb{R}^{1 \times \hat{n}}  & \mbox{ if }  D_l = 0
             \end{array}\right.,
  \end{split}
  \end{equation}
   where $C_s$ is the matrix formed by the first $n$ rows of $C_l$.
   Note that either $D_l$ is full column rank, and hence by positive definiteness of $Q,R$, $D_l^TSD_l$ is invertible, or $B_l=0,D_l=0,k=1$.  Hence, \eqref{DARE} is always
   well-defined.  For any $z \in \mathcal{V}(E,A,B)$, let $(x^*,u^*)$ be defined as
 \begin{equation}
 \label{dual:finhor:sol}
 \begin{split}
   & ({x^*}^T(s),{u^*}^T(s))^T=(C_{l}-D_{l}K(t-s))v(s)\,, \\
   & \dot v(s) = (A_l-B_lK(t-s))v(s) \mbox{ and \ \ } v(0)=\LMAP(z)\,.
 \end{split}
 \end{equation}
 Furthermore, define
  $K_f(s)=(C_u+D_uK(t-s))\LMAP E$, where $C_u$ and $D_u$ are the matrices formed by the last $m$ rows of $C_l$ and $D_l$ respectively. Define the matrices:
  \begin{equation}
 \label{finhor:feedback1}
     \begin{split}
     & K_1(t)=(C_l-D_lK(t_1-t))\LMAP E - \begin{bmatrix} I_n & 0 \\ 0 & 0 \end{bmatrix}\,, \\
     & K_2(t)=\begin{bmatrix} 0 & 0 \\ 0 & -I_m \end{bmatrix}.
    \end{split}
  \end{equation}
  \begin{Theorem}
  \label{finhopr:opt}
   With the notation above, $(x^*,u^*)$ is a solution
   of the finite horizon optimal control problem for the interval $[0,t]$ and the initial state $z$.
   The optimal value of the cost function is
  \begin{equation}
 \label{finhor:cost}
    J_t^{*}=J_{t}(x^*, u^*)=(\LMAP(z))^TP(t)(\LMAP(z)).
  \end{equation}
  Furthermore, $u^{*}(s)=K_f(s)x^*(s)$, $s \in [0,t]$, and
  $(x^{*},u^{*})=(x,u)$ is the unique (up to modification on a set of measure zero) solution of
  \begin{equation}
 \label{finhor:feedback2}
     \begin{split}
     & \dfrac{d(Ex(t))}{dt} =Ax(t)+Bu(t)\,, \\
     & K_1(t)x(t)+K_2(t)u(t) = 0,
    \end{split}
   \end{equation}
   on $[0,t_1]$, such that $Ex(0)=z$.\footnote{I.e., for any $(x,u) \in AC([0,t],\mathbb{R}^n) \times L^1_{loc}([0,t],\mathbb{R}^m)$ such that $Ex(0)=z$:
   $\dfrac{d(Ex(t))}{dt} =Ax(t)+Bu(t)$ a.e, and $K_1(t)x(t)+K_2(t)u(t) = 0$ a.e., if and only if $(x,u)=(x^*,u^*)$ a.e.}
 \end{Theorem}
 %Intuitively, Theorem \ref{finhopr:opt} relies on converting  Problem \ref{opt:fin} into a
 %finite horizon optimal control problem for the associated ODE-LTI.

 Next, we present the solution to the infinite horizon control problem. Just like in the classical case,
 we will need a certain notion of stabilizability for solvability of the infinite horizon LQ control problem. %The right notion is the following.
 \begin{Definition}[Behavioral stabilizability]
 \label{def:stabil}
  The DAE-LTI \eqref{dae:sys} is said to be behaviorally stabilizable from $z \in \mathcal{V}(E,A,B)$, if
  there exists $(x,u) \in \mathscr{D}_{z}(\infty)$ such that
  $\lim_{t \rightarrow \infty} x(t)=0$.
 \end{Definition}
%  We will show that terminology of \cite{BergerThesis,DAEBookChapter}, DAE \eqref{dae:sys} is said to be behaviorally stabilizable, if it is
%  behaviorally stabilizable for all $z \in \mathcal{V}(E,A,B)$.
% That is, the notion defined  in
% Definition \ref{def:stabil} is a generalization of behavioral stabilizability defined in
% \cite{BergerThesis,DAEBookChapter}. Similarly to behavioral stabilizability, it does not guarantee existence of a feedback control law which renders the closed-loop system stable. However, it guarantees existence of a stabilizing control law in the form of an algebraic constraint, see Remark \ref{interconnect}.
 Behavior stabilizability from $z$ can be interpreted in terms of the associated ODE-LTI as follows.
 Let
  $\mathscr{S}=(A_l,B_l,C_l,D_l)$ be an ODE-LTI associated
  with \eqref{dae:sys} and let $\LMAP$ be the corresponding state map. Let $\V_g$ denote the stabilizability subspace of
  $\mathscr{S}$.
 Recall from \cite{TrentelmanBook} that $\V_g$ is the set of all
 initial states $p_0$ of $\mathscr{S}$, for which there exists an input $g$
 such that the corresponding state trajectory $p$ starting from $p_0$
 has the property that $\lim_{t \rightarrow \infty} p(t)=0$.
 \begin{Lemma}
 \label{lemma:stabilizable}
  The DAE-LTI \eqref{dae:sys} is behaviorally
  stabilizable from $z \in \mathcal{V}(E,A,B)$ if and only if $\LMAP(z)$ belongs to the stabilizability subspace $\V_g$ of $\mathscr{S}$.
 \end{Lemma}
 In order to solve the infinite horizon control problem for DAE-LTIs, we reformulate it as
 an infinite horizon control problem for the associated ODE-LTIs. However, for ODE-LTIs, infinite horizon LQ control
 problems can be solved only for stabilizable ODE-LTIs. For this reason, we will need to define the restriction of
 an associated ODE-LTI to its stabilizability subspace.
 More precisely, consider the ODE-LTI $\mathscr{S}=(A_l,B_l,C_l,D_l)$ associated with \eqref{dae:sys} and consider
 its stabilizability subspace $\V_g$.
 From \cite{CalierDesoer} it then follows that $\V_g$ is $A_l$-invariant and $\IM B_l \subseteq \V_g$.
  Hence, there exists a basis transformation $T$ such that
  $T(\V_g)=\IM \begin{bmatrix} I_l & 0 \\ 0 & 0 \end{bmatrix}$,
  $l=\dim \V_g$ and in this new basis,
  \[
    \begin{split}
      & TA_lT^{-1}=\begin{bmatrix} A_g & \star \\ 0 & \star \end{bmatrix}, \quad
      TB_l=\begin{bmatrix} B_g \\ 0 \end{bmatrix}, \quad
       C_lT^{-1}=\begin{bmatrix} C_g^T \\ \star \end{bmatrix}^T, \quad
   \end{split}
  \]
  $A_g \in \mathbb{R}^{l \times l}, B_g \in \mathbb{R}^{l \times k},  C_g \in \mathbb{R}^{(n+m) \times l}$.
  Denote by $\mathscr{S}_g=(A_g,B_g,C_g,D_g)$,
  where $D_g=D_l$.
  \begin{Definition}[Stabilizable associated ODE-LTI]
  We call $\mathscr{S}_g$ a stabilizable ODE-LTI
  associated with \eqref{dae:sys} and we call $\LMAP_g=\begin{bmatrix} I_l & 0 \end{bmatrix} T\LMAP$ the associated state
  map.
 \end{Definition}
The ODE-LTI $\mathscr{S}_g$ represents the restriction of $\mathscr{S}$ to the subspace $\V_g$. It follows that $\mathscr{S}_g$ is stabilizable.
Moreover, since all associated ODE-LTIs of \eqref{dae:sys} are feedback equivalent, then all associated stabilizable ODE-LTIs of \eqref{dae:sys} are also feedback equivalent.
Consider a stabilizable ODE-LTI $\mathscr{S}_g=(A_g,B_g,C_g,D_g)$  associated with \eqref{dae:sys}, and the corresponding state map $\LMAP_g$ and
assume that $\LMAP(z) \in \V_g$.   Consider the following algebraic Riccati equation:
    \begin{equation}
    \label{ARE}
    \begin{split}
       & 0 = PA_g+A^T_gP- K^T(D^T_gSD_g)K + C^T_gSC_g\,, \\
       & K=\left\{ \begin{array}{rl}
          (D^T_gSD_g)^{-1}(B^T_gP+D^T_gSC_g)      & \mbox{ if } D_l \ne 0 \\
         0 \in \mathbb{R}^{1 \times l} & \mbox{ if } D_l=0.
       \end{array}\right. \\
     & S=\mathrm{diag}(Q,R)
    \end{split}
    \end{equation}
   Note that either $D_l$ is full column rank, and hence by positive definiteness of $Q,R$, $D_l^TSD_l$ is invertible, or $B_l=0,D_l=0,k=1$.  In the former case, $D_g=D_l$ is full column rank, in the
   latter case, $D_g=0$ and $B_g=0$
  Hence, \eqref{DARE} is always
   well-defined.
 \begin{Lemma}
 \label{opt:control:lemma1}
 The algebraic Riccati equation \eqref{ARE} has a unique symmetric  solution $P > 0$ and $A_g-B_gK$ is a stable matrix.
 \end{Lemma}
Consider now the tuple $(x^*,u^*) \in L^2_{loc}([0,+\infty),\mathbb{R}^n) \times L^2_{loc}([0,+\infty),\mathbb{R}^m)$ such that
 \begin{equation}
      \label{opt8}
        \begin{split}
        & \dot {v^*} = (A_g-B_gK)v^*\,, \quad v^*(0)=\LMAP_g(Ex_0)\,, \quad \\
        & ({x^*}^T,{u^*}^T)^T=(C_g-D_gK)v^*\,.
        \end{split}
 \end{equation}
  Furthermore, define the following matrices:
\begin{equation}
      \label{opt:feedback1}
  \begin{split}
  & K_f =(C_u-D_uK)\LMAP\,, \\
  & K_1=(C_g-D_gK)\LMAP E - \begin{bmatrix} I_n & 0 \\ 0 & 0 \end{bmatrix}\,, \\
  & K_2=\begin{bmatrix} 0 & 0 \\ 0 & -I_m \end{bmatrix}\,,
 \end{split}
\end{equation}
where $C_u$ and $D_u$ are the matrices formed by the last $m$ rows of $C_g$ and $D_g$ respectively. We can now state the following.
 \begin{Theorem}
 \label{opt:control}
 The following are equivalent:
 \begin{itemize}
 \item{\textbf{(i)}}
  The infinite horizon optimal control problem is solvable for $z \in \mathcal{V}(E,A,B)$
 \item{\textbf{(ii)}}
  The DAE-LTI \eqref{dae:sys} is behaviorally stabilizable from $z$
 %$\LMAP(z)$ belongs to the stabilizability subspace $\V_g$ of the associated ODE-LTI $\mathscr{S}$
 \item{\textbf{(iii)}}
  \( J^{*}_{\infty} = \limsup_{t \rightarrow \infty} \inf_{(x,u) \in \mathscr{D}_{z}(t)} J_t(x,u) < +\infty
  \)
\end{itemize}
 If either of  the conditions\textbf{(i)} -- \textbf{(ii)}  hold, then  $(x^*,u^*)$ from \eqref{opt8}
 is a solution of the infinite horizon optimal control problem for the initial state $z$. Moreover,
  \begin{equation}
      \label{opt9}
  J^{*}_{\infty}=J_{\infty}(x^*,u^*)=(\LMAP_g Ez)^T P\LMAP_g Ez\,,
  \end{equation}
  $u^*=K_fx^*$, $(x^{*},u^{*})$ is a solution of the DAE-LTI
  \begin{equation}
  \label{opt:feedback2}
     \begin{split}
     & \dfrac{d(Ex(t))}{dt} =Ax(t)+Bu(t)\,, \\
     & K_1x(t)+K_2(u(t) = 0
    \end{split}
  \end{equation}
   on $[0,+\infty)$ such that $Ex^*(0)=z$, and if $(x,u)$ is a solution of \eqref{opt:feedback2} on $[0,+\infty)$ such that $Ex(0)=z$, then $(x,u)=(x^*,u^*)$ a.e.
\end{Theorem}
 The proof of  Theorem \ref{opt:control} implies that
 in the formulation of optimal control problem, we can replace $\limsup$ by $\lim$.

  Note that the existence of solution for Problem \ref{opt:fin} and Problem \ref{opt:contr:def} and its computation
  depend only on the matrices $(E,A,B,Q,R,Q_0)$.
  Indeed, an ODE-LTI $\mathscr{S}$ associated with $(E,A,B)$
  can be computed from $(E,A,B)$, and the solution of the
  associated LQ problem can be computed using $\mathscr{S}$ and the
  matrices $Q,Q_0,R$.
  Notice that the only condition for the existence of a solution is
  behavioral stabilizability from $z$, and this can be checked by verifying if $\LMAP(z)$ belongs to the stabilizability
  subspace of $\mathscr{S}$.  The latter can be done by an algorithm.
  The Matlab code for solving Problem \ref{opt:fin} and Problem \ref{opt:contr:def} and checking behavioral stabilizability
  is available at \verb'http://sites.google.com/site/mihalypetreczky/'.

 We would like to conclude this section with a  short discussion on the notion of  stabilizability we proposed.
 First, there are several equivalent ways to define behavioral stabilizability. Below we state some of them.
%  From the proof of Lemma \ref{lemma:stabilizable} we can deduce the following equivalent formulations of behavioral stabilizability.
  \begin{Corollary}
  \label{col:stabilizable}
    For any $z \in \mathcal{V}(E,A,B)$ the following are equivalent.
  \begin{itemize}
  \item{(i)}   \eqref{dae:sys} is behaviorally stabilizable from $z$
  \item{(ii)}  there exist $(x,u) \in \mathscr{D}_{z}(\infty)$ such that $\lim_{t \rightarrow +\infty} Ex(t)=0$
%  \item {(iii)} there exists $(x,u) \in \mathscr{D}_{z}(\infty)$ such that $\lim_{t \rightarrow +\infty} x(t)=0$, $\lim_(t \rightarrow +\infty} u(t)=0$, and $x$ and $u$ are smooth functions.
  \item{(iii)} for any $(x,u) \in \mathcal{B}_{\mathbb{R}}(E,A,B)$ such that $Ex(0)=z$ and $x$ is absolutely continuous,
              there exist $(x_o,u_o) \in \mathcal{B}_{\mathbb{R}}(E,A,B)$ such that $(x_o(t),u_o(t))=(x(t),u(t))$ for all $t < 0$, $Ex(0)=z$ and
              $\lim_{t \rightarrow +\infty} x_o(t)=0$ ,$\lim_{t \rightarrow +\infty} u_o(t)=0$ and $x_o$ is absolutely continuous
  \end{itemize}
  \end{Corollary}
  In fact,   Part (iii) of Corollary \ref{col:stabilizable} implies that if \eqref{dae:sys} is behaviorally stabilizable for all $z \in \mathcal{V}(E,A,B)$, then \eqref{dae:sys} is behaviorally stabilizable in the sense of
  \cite{DAEBookChapter,BergerThesis}. Hence,
  \begin{Corollary}[{\cite[Corollary 4.3]{DAEBookChapter}},{\cite[Proposition 3.3]{BergerThesis}}]
  \label{col:stabilizable1}
   The DAE-LTI \eqref{dae:sys} is stabilizable for all $z \in \mathcal{V}(E,A,B)$ $\iff$ $\forall \lambda \in \mathbb{C}, \mathrm{Re} \lambda \ge 0: \mathrm{rank} \begin{bmatrix} \lambda E - A, & B \end{bmatrix} = nrank \begin{bmatrix} sE - A, & B \end{bmatrix}$.
   Here,  $nrank \begin{bmatrix} sE - A, & B \end{bmatrix}$ denotes the rank of the polynomial matrix $\begin{bmatrix} s E - A, & B \end{bmatrix}$ over the quotient field of polynomials in the variable $s$.
  \end{Corollary}
  Note that behavior stabilizability from all $z \in \mathcal{V}(E,A,B)$ is equivalent to behavioral
  stabilizability in the sense of \cite{DAEBookChapter}, and the latter is equivalent to existence of an
  algebraic constraint which stabilizes the closed-loop system, see \cite{DAEBookChapter}.
  By Theorem \ref{opt:control},
  behavior stabilizability from all $z \in \mathcal{V}(E,A,B)$ is equivalent to the existence of a
  solution of Problem \ref{opt:contr:def} for all $z \in \mathcal{V}(E,A,B)$.  The resulting optimal state trajectory
  $x^*$ converges to zero, and it can be enforced by adding the algebraic constraint $K_1x+K_2u=0$ to
  the original DAE-LTI. In fact,  $K_1x+K_2u=0$ from Theorem \ref{opt:control}
  is a particular instance of a stabilizing algebraic constraint from
  \cite{BergerThesis,DAEBookChapter}.
  However, as it was already pointed out in \cite{BergerThesis,DAEBookChapter},
  behavioral stabilizability does not imply existence of a stabilizing feedback. Below we present an
  example which is behaviorally stabilizable but cannot be stabilized by a state feedback.
%  We have shown that for the solution $(x^*,u^*)$ of the optimal control problem (both finite and infinite horizon), the optimal input $u^*$ can be expressed as a state-feedback.
%  Note, however, that the feedback law does not determine the control
%  input uniquely, since DAE~\eqref{dae:sys} may admit several solutions starting
%  from the same initial state.
  % If the DAE has at most one solution from any
  %initial state, in particular, if the DAE is regular, then the feedback law
  %above determines the optimal trajectory $x^{*}$ uniquely.
%  We argue that this is inherent to DAE-LTIs.
\begin{Example}
\label{interconnect}
    Consider
    the following DAE-LTI
 %$E=\begin{bmatrix} 1 & 0 & 0   \\ 0 & 1 & 0 \end{bmatrix}$, $A=\begin{bmatrix} 1 & 0 & 0 \\ 0 & 1 & 1 \end{bmatrix}$, $B=\begin{bmatrix} 1 & 0 \end{bmatrix}^T$
%Consider the following decomposition of state variables: $x=\begin{bmatrix} x_1 & x_2 & x_3 \end{bmatrix}$.
 %   Then this DAE-LTI can be rewritten as
    \begin{equation}
    \label{interconnect:example}
     \begin{split}
         & \dot x_1 = x_1 + u\,, \\
         & \dot x_2 = x_2 + x_3\,.
     \end{split}
    \end{equation}
  Since the second equation does not depend on $u$, no matter how we choose the feedback $u=g(x)$ for some function $g$,
  it will not influence $x_2$. That is, there is no chance to enforce any restriction on $x_2$by using the control input
  only.
 %the solution of the system with this control law will admit several solutions from any initial state. In fact, it will always admit a solution $(x,u)$ for which $J_{\infty}(x,u)$ is infinite.
  However, optimal and stabilizing control is still possible.
Consider the matrices $Q=I_3$, $R=1$, $Q_0=I_2$ and the corresponding optimal control problem (Problem \ref{opt:contr:def}). An associated ODE-LTI  $(A_l,B_l,C_l,D_l)$ can be chosen so that
  $A_l=I_2,B_l=\begin{bmatrix} 0 & 1 \\ 1 & 0 \end{bmatrix}$, $C_l=\begin{bmatrix} I_2 \\ 0 \end{bmatrix}$, $D_l=\begin{bmatrix} 0 \\ I_2 \end{bmatrix}$. This system is clearly controllable and hence stabilizable, and
  thus it can be taken as a stabilizable ODE-LTI associated with the DAE-LTI. The solution of the Riccati equation and the corresponding matrices $K_1$ and $K_2$ can readily be computed.
 \end{Example}
  The optimal control problem
  described in Example \ref{interconnect} arises when trying to solve the  problem of estimating the state $z_1$
  of the following noisy DAE with the output $y$: $\dot z_1 = -z_1+f_1$, $\dot z_2 = -z_2+f_2$, $z_2=f_3$, $y=z_1+\eta$.
  Here $f_1,f_2,f_3,\eta \in L^2([0,+\infty),\mathbb{R})$ are deterministic noise signals such that: $z^2_1(0)+\int_0^{\infty} f_1^2(s)+f_2^2(s)+f_3^2(s)+\eta^2(s)ds \le 1$, i.e. the energy of
  $f_1,f_2,f_3,\eta$ is bounded and the unknown initial state $z_1(0)$ is bounded. Then
  according to \cite{ZhukPetrezky2013TAC,ZhukPetreczkyCDC13}, in order to construct an estimate $\hat z_1$ of $z_1$ from $y$ with the minimal worst-case estimation error, one needs to solve the LQ problem of Example \ref{interconnect}. Conversely, if $\hat z_1$ exists, then the LQ control problem from Example \ref{interconnect} will have a solution.
  That is, even for such toy models (if $\eta=0$, the problem is trivial), the state estimation problem yields an optimal control problem which cannot be solved by state
  feedback alone.

  Example \ref{interconnect} shows that DAE-LTIs which are behaviorally stabilizable but do not admit a stabilizing feedback occur naturally.
  It shows that the lack of a stabilizing feedback control is not a shortcoming of the definition,
  but a sign that concept of the feedback control might be too restrictive for DAE-LTIs. In our opinion, one should consider more general controllers, for example, controllers which are represented by algebraic constraints $K_1x+K_2u=0$. The latter can be viewed as a controller, if we follow the philosophy of J.C. Willems
  \cite{WillemsInterconnect,WillemsMagazine}.
   Note that many physical control devices cannot be described as feedback controllers, \cite{WillemsInterconnect,WillemsMagazine}, including such simple example as mass-spring-dumper systems and
  electrical circuits. Note that a classical feedback $u=Kx$ is just a specific case of a controller enforcing algebraic constraints: it can be represented as $u-Kx=0$.

  The fact that we consider DAE-LTIs which cannot be optimized or even stabilized by state feedback explains
  why our results differ from \cite{Kurina1,Kurina2,Kurina3,Kurina4}.
  In  \cite{Kurina1,Kurina2,Kurina3,Kurina4} sufficient conditions for existence of an optimal state feedback control were presented. In particular, the conditions of \cite{Kurina1,Kurina2,Kurina3,Kurina4} imply existence of a
 stabilizing state feedback control law. These conditions cannot be satisfied by systems which cannot be stabilized by
  state feedback alone.  The system from Example \ref{interconnect} is one such system,
  and for that system the conditions of \cite{Kurina1,Kurina2,Kurina3,Kurina4} never hold, no matter which quadratic cost function we choose.
  %Note that such optimal control problems arise naturally when trying to solve minimax observer design problems, \cite{ZhukPetreczkyCDC13,ZhukPetrezky2013TAC}.
 % \end{Remark}

\subsection{Proofs}
\label{lqr:control:proof}
 In order to present the proofs of Theorem \ref{finhopr:opt} and Theorem \ref{opt:control}, we rewrite Problems \ref{opt:fin} -- \ref{opt:contr:def} as LQ control problems for ODE-LTIs.
 To this end, consider an ODE-LTI $\mathscr{S}=(A_l,B_l,C_l,D_l)$
 and let $\LMAP=(EC_s)^{+}$ be the corresponding state map. Recall that $C_s$ is the matrix formed by the first $n$ rows of $C_l$.
  Consider the following linear quadratic control problem. For every initial state $v_0$, for every interval $I$ containing $[0,t_1]$ and for every
  $g \in L^1_{loc}(I,\mathbb{R}^k)$ define the cost functional $\mathscr{J}_t(v_0,g)$:
  \begin{equation}
  \label{lpg:eq:cost1}
    \begin{split}
     & \mathscr{J}_{t}(v_0,g) =  v^T(t)C_s^TE^T Q_0 EC_{s} v(t) + \\
     & +\int_0^{t} \nu^T(s)\begin{bmatrix} Q & 0 \\ 0 & R \end{bmatrix} \nu(s)ds\,,    \\
     & \dot v = A_lv + B_l g\,,  \mbox{ and }  v(0)=v_0\,, \\
     & \nu = C_lv+D_lg\,.
    \end{split}
   \end{equation}
 For  any $g \in L^1_{loc}([0,+\infty),\mathbb{R}^k)$ and $v_0 \in \mathbb{R}^{\hat{n}}$, define
 \begin{equation}
 \label{lpg:eq:cost2}
    \mathscr{J}_{\infty}(v_0,g)=\limsup_{t \rightarrow \infty} \mathscr{J}_{t}(v_0,g).
  \end{equation}
  In our next theorem we prove that the problem of minimizing the cost function
  $J_{t}$ for \eqref{dae:sys} is equivalent to minimizing the cost function $\mathscr{J}_{t}$ for
  the associated ODE-LTI $\mathscr{S}$.
  \begin{Theorem}
 \label{lqg:eq}
    With the notation above, let $z \in \mathcal{V}(E,A,B)$, $I=[0,t]$, $t > 0$ or $I= [0,+\infty)$. For $g \in L^1_{loc}(I,\mathbb{R}^k)$ denote by $\nu_{\mathscr{S}}(v_0,g)$ the output trajectory of
  the associated ODE-LTI $\mathscr{S}$, which corresponds to the initial state $v_0$  and input $g$.

 \textbf{(i)}
     For any $g \in L^1_{loc}(I,\mathbb{R}^{k})$, and
     for any $(x,u) \in \mathscr{D}_{z}(t)$, such that
     $(x^T,u^T)^T=\nu_{\mathscr{S}}(\LMAP(z),g)$ a.e. ,
     \begin{equation}
      \label{daelincost}
     \begin{split}
        & J_{t}(x,u)=\mathscr{J}_t(\LMAP(z),g)\,, \\
         & J_{\infty}(x,u)=\mathscr{J}_{\infty}(\LMAP(z),g)\,.
     \end{split}
     \end{equation}
 \textbf{(ii)}
   For all $t  \in (0,+\infty)$, $(x^*,u^*) \in \mathscr{D}_{z}(t)$ is a solution
   of the finite horizon optimal control problem if and only if
   there exists $g^* \in L^1([0,t],\mathbb{R}^k)$ such that
   $({x^*}^T,{u^*}^T)^T=\nu_{\mathscr{S}}(\LMAP(z),g^*)$ a.e. and
   \begin{equation}
     \label{eq:LQfin}
       \mathscr{J}_{t}(v_0,g^*)=\inf_{g \in L^1_{loc}([0,t],\mathbb{R}^k)} \mathscr{J}_{t}(v_0, g) < +\infty\,.
   \end{equation}
 \textbf{(iii)}
  The tuple $(x^{*},u^{*}) \in \mathscr{D}_{x_0}(\infty)$ is a solution of
  the infinite horizon optimal control problem if and only if there
  exists an input $g^* \in L^1_{loc}([0,+\infty), \mathbb{R}^k)$ such that
  $({x^*}^T,{u^*}^T)^T=\nu_{\mathscr{S}}(\LMAP(z),g^*)$ a.e., and
  \begin{equation}
  \label{eq:LQinf}
     \mathscr{J}_{\infty}(v_0,g^{*}) = \limsup_{t \rightarrow \infty} \inf_{g \in L^1_{loc}([0,t_1],\mathbb{R}^k)} \mathscr{J}_{t}(v_0,g) < +\infty.
  \end{equation}
 \end{Theorem}
 \begin{IEEEproof}[Proof of Theorem \ref{lqg:eq}]
  Equation \eqref{daelincost} follows by routine manipulations and by
  noticing that the first $n$ rows of $C_lv(t)+D_lg(t)$ equal
  $C_sv(t)+D_sg(t)$ and as $ED_s=0$ (see Definition \ref{def:asscoLTI} for the definitions of $C_s$ and $D_s$),
  $E(C_sv(t)+D_sg(t))=EC_sv(t)$.
  The rest of theorem follows by noticing that for any element $(x,u)$
  $\mathscr{D}_{z}(t_1)$ or $\mathscr{D}_{z}(\infty)$, there exist
  $g \in L^1_{loc}(I,\mathbb{R}^k)$ such that for $v=\LMAP Ex$, $(x^T,u^T)^T=C_lv+D_lg$ a.e,  $\dot v = A_lv+B_lv$, a.e.  $v(0)=\LMAP Ex(0)=\LMAP z$, and conversely,
  for any $(v,g) \in AC(I,\mathbb{R}^{\hat{n}}) \times L^1_{loc}(I,\mathbb{R}^k)$ such that $\dot v = A_lv+B_lg$, $v(0)=\LMAP z$, $C_lv+D_lg \in \mathscr{D}_z(t_1)$ (if $I=[0,t_1]$) or
  $C_lv+D_lg \in \mathscr{D}_{z}(\infty)$ (if $I=[0,+\infty)$).
 \end{IEEEproof}
 The proof of Theorem \ref{finhopr:opt} can then be derived from the classical results (see \cite{KwakernaakBook}).
  \begin{IEEEproof}[Proof of Theorem \ref{finhopr:opt}]
   By definition of an associated ODE-LTI, there are two possible cases:
   either $D_l$ is full column rank or $D_l=0$. This two cases cover all the possibilities. We start with the case when $D_l \ne 0$ is full column rank.

   Let us first apply the feedback transformation $g=\hat{F}v+Uw$ to $\mathscr{S}=(A_l,B_l,C_l,D_l)$
  with $U=(D^T_lSD_l)^{-1/2}$ and $\hat{F}=-(D^T_lSD_l)^{-1}D^T_lSC_l$,
  as described in \cite[Section 10.5, eq. (10.32)]{TrentelmanBook}.
 Note that $D_l$ is injective and hence $U$ is well defined.
  Consider the linear system
  \begin{equation}
  \label{opt:st1}
   \begin{split}
   & \dot v = (A_l+B_l\hat{F})v+B_lUw \mbox{ and } v(0)=v_0 \,.
   \end{split}
  \end{equation}
  For any $w \in L^1_{loc}([0,t],\mathbb{R}^k)$, the state trajectory $v$ of \eqref{opt:st1} equals the state trajectory of
  $\mathscr{S}=(A_l,B_l,C_l,D_l)$ for the input $g=\hat{F}v+Uw$ and initial state $v_0$. Moreover,
from Theorem~\ref{theo:geo_contr} it follows that
all inputs $g \in L^1_{loc}([0,t],\mathbb{R}^k)$ of $(A_l,B_l,C_l,D_l)$ can be represented in such a way. Define
   \begin{equation*}
  %\label{opt:st2}
  \begin{split}
      & \widehat{\mathscr{I}_t}(v_0,w)= v^T(t){C}_s^TE^TQ_0 EC_{s} v(t)+ \\
       & +\int_0^t (v^T(s)(C_l+D_l\hat{F})^TS(C_l+D_l\hat{F})v(s)+w^T(s)w(s))ds,
   \end{split}
  \end{equation*}
  where $v$ is a solution of \eqref{opt:st1}, and $C_s$ is the matrix formed by the first $n$ rows of $C_l$.
  It is easy to see that $\mathscr{J}_t(v_0,g)=\widehat{\mathscr{I}}_t(v_0,w)$ for $g=\hat{F}v+Uw$ and any initial state $v_0$ of $\mathscr{S}$.

  Consider now the problem of minimizing $\widehat{\mathscr{I}}_t(v_0,w)$. The solution of this problem can be found using \cite[Theorem 3.7]{KwakernaakBook}.  Notice that \eqref{DARE} is equivalent to the Riccati differential equation described in \cite[Theorem 3.7]{KwakernaakBook} for the problem of minimizing $\widehat{\mathscr{I}}_t(v_0,w)$. Hence, by \cite[Theorem 3.7]{KwakernaakBook}, \eqref{DARE} has a unique positive solution $P$, and for the optimal input $w^*$, $g^{*}=\hat{F}v^{*}+Uw^{*}=-K(t_1-t)v(t)$ satisfies~\eqref{eq:LQfin}, and \( \dot v(t) = (A_l-B_lK(t_1-t))v(t) \mbox{ and \ \ } v(0)=v_0 \). From Theorem \ref{lqg:eq} and Definition \ref{def:asscoLTI} it then follows that $({x^*}^T,{u^*}^T)^T=C_lv^*+D_lg^*$ is the solution of the Problem 1 and that \eqref{finhor:cost} holds.

  Assume now that $D_l=0$. From the definition of an associated ODE-LTI is then follows that $B_l=0$. Hence, the associated ODE-LTI is in fact an autonomous system. Define now
  $P(t)= e^{A_l^Tt}(EC_s)^TQ_0EC_se^{A_lt} + \int_0^t e^{A^T_ls}C_l^TSC_le^{A_ls}ds$. It is then easy to see that $P(t)$ satisfies \eqref{DARE}.
  %Moreover, since $D_l=0$, the set $\mathcal{D}_{z}(t)$ contains exactly one element $(x^*,u^*)$ and this tuple satisfies \eqref{dual:finhor:sol}.
  %Note that since $(x^*,u^*)$ is the only element of $\mathcal{D}_{z}(t)$, $J^*_t=J_t(x^*,u^*)$ and thus $(x^*,u^*)$ is necessarily a solution of the finite horizon optimal control problem for the interval $[0,t]$ and initial state $z$. Finally, notice that the state trajectory $v$ from \eqref{dual:finhor:sol} satisfies
  Let $(x^*,u^*) \in \mathscr{D}_z(t)$  be such that \eqref{dual:finhor:sol} is satisfied. Since $D_l=0,B_l=0$,
  for any $(x,u) \in \mathscr{D}_z(t)$, $x=x^*$ a.e. and $u=u^*$ a.e. Moreover, since $Ex$ and $Ex^*$ are absolutely continuous,
  $Ex=Ex^*$ a.e. implies $Ex=Ex^*$. Hence, $J_t(x,u)=J_t(x^*,u^*)$ and thus
  $(x^*,u^*)$ is necessarily a solution of the finite horizon optimal control problem for the interval $[0,t]$ and initial state $z$.
  Finally, notice that the state trajectory $v$ from \eqref{dual:finhor:sol} satisfies
  $v(s)=e^{A_ls}\LMAP(z)$ and $({x^{*}}^T,{u^*}^T)^T(s)=C_le^{A_ls}\LMAP z$ and hence \eqref{finhor:cost} holds.

  Finally, in both cases ($D_l$ is full rank or $D_l=0$), $u^{*}(s)=K_f(s)x^*(s)$ for all $s \in [0,t]$ and all the solutions $(x,u)$ of \eqref{finhor:feedback2}  on $[0,t]$ are those which satisfy
  $(x^T,u^T)^T(s)=(C_{l}-D_{l}K(t-s))v(s)$ a.e., $\dot v(s) = (A_l-B_lK(t-s))v(s)$, $v(s)=\LMAP(Ex(s))$, $s \in [0,t]$. Hence, $(x^*,u^*)$ is indeed the only solution of \eqref{finhor:feedback2} such that
  $Ex^*(0)=z$.
 \end{IEEEproof}
  %Next, we present the proof of Lemma \ref{lemma:stabilizable} and Corollary \ref{col:stabilizable}.
  \begin{IEEEproof}[Proof Lemma \ref{lemma:stabilizable}]
  \textbf{"if part"}
  Assume that $\mathscr{S}=(A_l,B_l,C_l,D_l)$ is
  stabilizable from $\LMAP(z)$.  Let $\V_g$ be the stabilizability subspace
  of $\mathscr{S}$.
  It then follows from \cite{TrentelmanBook}
  that $\LMAP(z) \in \V_g$ and
  there exists a feedback $F_l$ such that the restriction of
  $A_l+B_lF_l$ to $\V_g$ is stable and hence for any $v_0=\LMAP(z)$, there
  exists $v \in AC([0,+\infty),\mathbb{R}^{\hat{n}})$, such that  $\dot v = (A_l+B_lF_l)v$, $v(0)=v_0$,  and
  $\lim_{t \rightarrow \infty} v(t)=0$. Consider the output $(x^T,u^T)^T=(C_l+D_lF_l)v$.
  It then follows that $Ex(0)=z$ and $(x,u)$ is a solution of the DAE-LTI.
  Moreover, $\lim_{t \rightarrow \infty} x(t)=\lim_{t \rightarrow \infty} (C_{s}+D_{s}F_l)v(t)=0$. That is, DAE-LTI is stabilizable from $z$.

  \textbf{"only if part"}
 % Choose an asscociated ODE-LTI $(A_l,B_l,C_l,D_l)$.
  Assume that the DAE-LTI is stabilizable from $z$, and let $(x,u) \in \mathscr{D}_{z}(\infty)$ be such that $\lim_{t \rightarrow \infty} x(t)=0$. It then follows that there exist an input
  $g \in L^1_{loc}([0,+\infty),\mathbb{R}^k)$ such that
  $\dot v =A_lv + B_lg$, $(x^T,u^T)^T=C_lv+D_lg$ a.e., $v(0)=v_0=\LMAP(Ex_0)$ and
  $\LMAP(Ex(t))=v(t)$, $t \in [0,+\infty)$.
  In particular, $\lim_{t \rightarrow \infty} v(t)=\LMAP E \lim_{t \rightarrow \infty} x(t)=0$.
  That is, there exists an input $g$, such that the corresponding
  state trajectory $v$ of $\mathscr{S}$ starting from $\LMAP(z)$ converges
  to zero.
%  That is, for any initial state $v_0$ of $(A_l,B_l,C_l,D_l)$ there exists an input $g$ such that the corresponding state $v$ is stable.
  But this is precisely the definition of stabilizability of $(A_l,B_l,C_l,D_l)$
  from $\LMAP(z)$.
 \end{IEEEproof}
   \begin{IEEEproof}[Proof of Corollary \ref{col:stabilizable}]
   The implications (i) $\implies$ (ii) is trivial The implication (ii) $\implies$ (i) follows by noticing that in the proof of the "only if" part of Lemma \ref{lemma:stabilizable} it is
   sufficient to assume that $(x,u) \in \mathscr{D}_{z}(\infty)$ is such that $\lim_{t \rightarrow +\infty} Ex(t)=0$.

   (i) $\implies$ (iii) can be shown as follows.
   Let $\mathscr{S}$ be an associated ODE-LTI of \eqref{dae:sys}, $\LMAP$ be the corresponding state map and let $\V_g$ be the stabilizability subspace of $\mathscr{S}$. If \eqref{dae:sys} is stabilizable from
   $z$, then by Lemma \ref{lemma:stabilizable}, $\LMAP(z) \in \V_g$ and there exists a feedback control law $F_l$ such that $(A_l+B_lF_l)\V_g \subseteq \V_g$ and the restriction of $(A_l+B_lF_l)$ to $\V_g$ is stable. If $(x,u) \in \mathcal{B}_{\mathbb{R}}(E,A,B)$, then $\LMAP(Ex)=p \in AC(I,\mathbb{R}^{\hat{n}})$ is a solution of $\dot p = A_lp+B_lg$ for
   some $g  \in L^1_{loc}(\mathbb{R},\mathbb{R}^k)$ and
   $(x^T,u^T)^T=C_lp+D_lg$ a.e..  Let $C_s$, $D_s$ be the matrices formed by the first $n$ rows of $C_l$, $D_l$.
   If $x$ is absolutely continuous, then $x=C_sp+D_sg$ a.e. implies $x-C_sp=D_sg$ a.e. and $x-C_sp$ is absolutely continuous. Hence, by modifying
   $g$ on a set of measure zero, without loss of generality, we can assume that $D_sg$ is absolutely continuous. Let $\phi:\mathbb{R} \rightarrow [0,1]$ be a smooth function such that $\phi(t)=1$ for $t \le 0$ and $\phi(t)=0$ for $t \ge 1$.
   The existence of such a function follows from partition of unity, \cite{Diff:Geo}.
   Let $p_o \in AC(\mathbb{R},\mathbb{R}^{\hat{n}})$ be the solution of $\dot p_o=(A_l+B_lF_l)p_o+\phi B_l(g-F_l p_o)$ a.e, $p_o(0)=\LMAP(z)$.
   Notice that $p_o(t) \in \V_g$ for all $t \in \mathbb{R}$, since the stabilizability subspace of any ODE-LTI is invariant under the dynamics of this
   ODE-LTI. Notice that for all $t \le 0$,
   $\dot p_o(t)=A_lp_o(t)+B_lg(t)$ and hence $p_o|_{(-\infty,0]}=p|_{(-\infty,0]}$ by uniqueness of solutions of differential equations.
   Note that for $t \ge 1$, $\dot p_o(t)=(A_l+B_lF_l)p_o(t)$ and hence $\lim_{t \rightarrow +\infty} p_o(t)=0$.
   Set $g_o(t)=\phi(t)(g(t)-F_lp_o(t))+F_lp_o(t)$, $t \in \mathbb{R}$. It is clear that $g_o|_{(-\infty,0]}=g|_{(-\infty,0]}$, $D_sg_o$ is absolutely continuous and
   $\dot p_o=A_lp_o+B_lg_o$ a.e.
   Define $(x_o^T,u_o^T)^T = C_lp_o+D_lg^*$. It then follows that $(x_o,u_o) \in  \mathcal{B}_{\mathbb{R}}(E,A,B)$ and $Ex_o(0) = \LMAP EC_s  p(0)=z$,
   $(x_o(t),u_o(t))=(x(t),u(t))$ for all $t < 0$ and $x_o$ is absolutely  continuous.
   Moreover, $(x_o^T(t),u_o^T(t))^T=(C_l+D_lF_l)p_o(t)$ for all $t \ge 1$ and hence $x_o(t),u_o(t)$ converge to zero as $t \rightarrow +\infty$.

   The implication (iii) $\implies$ (i) can be shown as follows. Since $z$ is differentiably consistent, from Corollary \ref{col3:important} it follows that
    there exist $(x,u) \in \mathcal{B}_{\mathbb{R}}(E,A,B)$ such that $Ex(0)=z$ and $x$ is differentiable.
    Then from (iii)  it follows that there exist $(x_o,u_o) \in \mathcal{B}_{\mathbb{R}}(E,A,B)$ such that $(x_o(t),u_o(t))=(x(t),u(t))$ for all $t < 0$, $Ex(0)=z$ and
              $\lim_{t \rightarrow} (x^T_o(t),u^T_o(t))^T=0$. In particular, $(x_o|_{[0,+\infty)}, u_o|_{[0,+\infty)}) \in \mathscr{D}_{z}(\infty)$ and $\lim_{t \rightarrow +\infty} x_o(t)=0$, i.e.\eqref{dae:sys}
     is behaviorally stabilizable from $z£$.
  \end{IEEEproof}
 %   Notice that the set $\V_g$ is invariant with respect to the dynamics of $(A_l,B_l,C_l,D_l)$, i.e.
 % any state trajectory starting in $\V_g$ stays in $\V_g$. If then follows that if $\LMAP(z) \in \V_g$ and $(x,u) \in \mathscr{D}_z(I)$ for some interval $I$, then
 % $p(t)=\LMAP(Ex(t)) \in \V_g$, where $p$ and $(x,u)$ satisfy $\dot p = A_lp+B_lg$ a.e. $(x^T,u^T)^T=C_lp+D_lg$ for some  $g \in L_{loc}^1(I,\mathbb{R}^k)$.
 % That is, if \eqref{dae:sys} is stabilizable from $z$, it is stabilizable from any state which can be reached from $z$.
 %Notice that if $B_l=0$, then $\V_g$ is the stability subspace of $A_l$, i.e. it is the largest $A_l$ invariant subspace of $\mathbb{R}^{\hat{n}}$ such that the restriction of $A_l$ to
 % $\V_g$ has all eigenvalues in the left-half plane.
%Using Theorem \ref{lqg:eq} and the classical results on infinite horizon LQ control, we can solve the infinite horizon optimal control problem as follows.
 %Finally, we present the proof of Lemma \ref{opt:control:lemma1} and Theorem \ref{opt:control}.
 \begin{IEEEproof}[Proof of Lemma \ref{opt:control:lemma1}]
  If $D_g=0$, then $B_g=0$ and as $(A_g,B_g)$ is stabilizable, $A_g$ is stable. Then existence of $P=P^T > 0$ satisfying \eqref{ARE} follows from the existence of the observability grammian for a stable linear systems.

  Assume now that $D_g$ is full column rank.
  Let us apply the feedback transformation $g=\hat{F}v+Uw$ to $\mathscr{S}_g=(A_g,B_g,C_g,D_g)$
  with $U=(D^T_gSD_g)^{-1/2}$ and $\hat{F}=-(D^T_gSD_g)^{-1}D^T_gSC_g$, as
  described in \cite[Section 10.5, eq. (10.32)]{TrentelmanBook}.
To this end, notice that $(A_g+B_g\hat{F},B_gU)$ is stabilizable and $(S^{1/2}(C_g+D_g\hat{F}),A_g+B_g\hat{F})$ is observable. Indeed, it is easy to see that stabilizability of $(A_g,B_g)$ implies that of $(A_g+B_g\hat{F},B_gU)$.
   Observability of $(S^{1/2}(C_g+D_g\hat{F}),A_g+B_g\hat{F})$ can be derived as follows. Recall from Definition \ref{def:asscoLTI} that $EC_s$ is of full column rank
   and $ED_s=0$. Note that $D_g=D_l$. Let $\hat{C}_s$ be the matrix formed by the first $n$ rows of $C_g$. Then $E\hat{C}_s$ is the the restriction of the map $E(C_s+D_s\hat{F})=EC_s$
 to $\V_g$, hence $E\hat{C}_s$ is of full column rank, if $EC_s$ is injective. The latter is the case according to Definition \ref{def:asscoLTI}.
   Hence, $E(\hat{C}_s+\hat{D}_s\hat{F})=E\hat{C}_s$ is of full column rank, and
   thus the pair $(\hat{C}_s+\hat{D}_s\hat{F},A_g+B_g\hat F)$ is observable.

  Consider the ODE-LTI
  \begin{equation}
  \label{opt:st2}
   \begin{split}
   & \dot v = (A_g+B_g\hat{F})v+B_gUw \mbox{ and } v(0)=\hat{v}_0\,. \\
   \end{split}
  \end{equation}
  For any $w \in L^1_{loc}(I,\mathbb{R}^k)$, where $I=[0,t]$, $0 < t \in \mathbb{R}$ or $I=[0,+\infty)$,
  the state trajectory $v$ of \eqref{opt:st2} equals the state trajectory of
  $\mathscr{S}_g=(A_g,B_g,C_g,D_g)$ for the input $g=\hat{F}v+Uw$ and initial state $\hat{v}_0$.
  Moreover, all inputs $g \in L^1_{loc}(I,\mathbb{R}^k)$ of $(A_g,B_g,C_g,D_g)$ can be represented in such a way.
  Define now
   \begin{equation*}
  %\label{opt:st2}
  \begin{split}
      & \widehat{\mathscr{I}}_t(\hat{v}_0,w)= v^T(t)(E\hat{C}_s)^TQ_0 E\hat{C}_{s} v(t)+ \int_0^t [w(s)^Tw(s) + \\
     &  +(v^T(s)(C_g+D_g\hat{F})^TS(C_g+D_g\hat{F})v(s)]ds,
   \end{split}
  \end{equation*}
  where $v$ is a solution of \eqref{opt:st2}.
   Consider now the problem of minimizing $\lim_{t \rightarrow \infty} \widehat{\mathscr{I}}_t(\hat{v}_0,w)$.
   Notice that \eqref{ARE} is equivalent to the algebraic Riccati equation
   described in \cite[Theorem 3.7]{KwakernaakBook} for the ODE-LTI \eqref{opt:st2} and for the
   infinite horizon cost function $\lim_{t \rightarrow \infty} \widehat{\mathscr{I}}_t(\hat{v}_0,w)$.
 Hence, by \cite[Theorem 3.7]{KwakernaakBook},
 \eqref{ARE}  has a unique positive definite solution $P$,
  and $A_g+B_g\hat{F}-B_gUU^TB_g^TP=A_g-B_gK$ is a stable matrix.
 \end{IEEEproof}
 \begin{IEEEproof}[Proof of Theorem \ref{opt:control}]
\textbf{(i) $\implies$ (ii)}
  If $(x^*,u^*)$ is a solution of the infinite horizon optimal control problem, then by Theorem \ref{lqg:eq}, there exists an input $g^* \in L^1_{loc}([0,+\infty),\mathbb{R}^k)$
such that $\mathscr{J}_{\infty}(\LMAP(z),g^*) <+\infty$. Let $v_0=\LMAP(z)$.
We claim that if $\mathscr{J}(v_0,g^*) < +\infty$,
then $\lim_{t \rightarrow \infty} v^*(t)=0$ for the state trajectory $v^*$ of
$\mathscr{S}$ which corresponds to the input $g^*$ and starts from $v_0$.
The latter is equivalent to $v^*(0)=v_0=\LMAP(Ex) \in \V_g$.
Let us prove that $\lim_{t \rightarrow \infty} v^*(t)=0$.
To this end, notice that $\mathscr{J}_{\infty}(v_0,g^*) < +\infty$ implies
 \( \int_0^{\infty} \nu(t)^T \begin{bmatrix} Q & 0 \\ 0 & R \end{bmatrix} \nu(t)dt < +\infty\), and as $\begin{bmatrix} Q & 0 \\ 0 & R \end{bmatrix}$
is positive definite, it follows that
$\mu\int_0^{\infty} \nu^T(t)\nu(t)dt < \int_0^{\infty} \nu(t)^T \begin{bmatrix} Q & 0 \\ 0 & R \end{bmatrix} \nu(t)dt < +\infty \)
for some $\mu > 0$ and so $\nu \in L^2([0,+\infty),\mathbb{R}^{n+m})$.
Consider the decomposition $\nu(t)=(x^T(t),u^T(t))^T$,
where $x(t) \in \mathbb{R}^n$.
By Theorem~\ref{dae2lin:theo_main} it follows that $\LMAP(Ex(t))=v^*(t)$ and
hence $v^* \in L^2([0,+\infty),\mathbb{R}^{\hat{n}})$. As $g^*$ is a linear function of $x,u$ and $v$
(see~Theorem~\ref{dae2lin:theo_main}) it follows that $g^* \in L^2([0,+\infty),\mathbb{R}^k)$.
Recalling that $\dot v^*(t)= A_lv^*(t)+B_lg^*(t)$ we write:
$v^*(t)=v^*(\tau)+A_l\int_\tau^t v^*(s)ds+B_l\int_\tau^t g^*(s)ds$ for $\tau<t$.
Since $\int_\tau^t \|v^*(s)\|ds\le(t-\tau)^\frac 12 \|v^*\|_{L^2([0,+\infty),\mathbb{R}^{\hat{n}})}$ where $\|v^*\|^2_{L^2([0,+\infty),\mathbb{R}^{\hat{n}})}:=\int_0^{\infty}\|v^*(s)\|_{\mathbb R^{\hat{n}}}^2 ds$,
 it follows that: $\|v^*(t)-v^*(\tau)\|_{\mathbb R^{\hat{n}}}\le (t-\tau)^\frac 12 (\|A_l\|\|v\|_{L^2([0,+\infty),\mathbb{R}^{\hat{n}})}+\|B_l\|\|g^*\|_{L^2([0,+\infty),\mathbb{R}^k)})$.
Hence $v^*$ is uniformly continuous. This and $v^* \in L^2([0,+\infty),\mathbb{R}^{\hat{n}})$ together with Barbalat's lemma imply that $\lim_{t \rightarrow \infty} v^*(t) = 0$. Hence the optimal state trajectory $v^*$ converges to zero, i.e. $g^*$ is a stabilizing control for $v_0$, and thus $v_0=\LMAP(z) \in \V_g$.

\textbf{(ii) $\implies$ (i)}
  Assume now that $v_0=\LMAP(z) \in \V_g$. Let us recall the construction of the ODE-LTI $\mathscr{S}_g$. It then follows that
  $T(\V_g)=\IM \begin{bmatrix} I_l & 0 \\ 0 & 0 \end{bmatrix}$, where
  $l=\dim \V_g$. Define the map
  $\Pi=\begin{bmatrix} I_l & 0 \end{bmatrix}T$. It then follows that
  for any $v_0 \in V_g$, and for any $g \in L^1_{loc}([0,+\infty],\mathbb{R}^k)$, $(x^T,u^T)^T$ is the output of $\mathscr{S}$ and $t\mapsto v(t)$ is the state trajectory of $\mathscr{S}$ starting from $v_0$ and driven by the input $g$, if and only if $(x^T,u^T)^T$ is the output of $\mathscr{S}_g$ and $t\mapsto\Pi(v)(t)$ is the state trajectory of $\mathscr{S}_g$ starting from $\Pi(v_0)$ and driven by $g$.
%In other words, the state trajectories of $\mathscr{S}_g$ correspond to the state trajectories of $\mathscr{S}$
%  which start in $\V_g$, and for these trajectories $\mathscr{S}_g$ and $\mathscr{S}$ produce the same output.
  For any initial state $\hat{v}_0$ of
  $\mathscr{S}_g$ define now the cost function
  \( \mathscr{I}_t(\hat{v}_0,g) \) as
  \[
    \begin{split}
      \mathscr{I}_t(\hat{v}_0,g) &=  v^T(t)(E\hat{C}_s)^T Q_0 E\hat{C}_{s} v(t) +\int_0^{t} \nu^T(s)S\nu(s)ds\,,    \\
     \dot v = A_gv + &B_g g\,,  \mbox{ and }  v(0)=\hat{v}_0\,, \\
     \nu = C_gv+&D_gg , \quad S=\mathrm{diag}(Q, R)\,,
    \end{split}
   \]
  where $\hat{C}_s$ is the matrix formed by the first $n$ rows of $C_g$.
 Notice that from the definition of $C_g$ it follows that
 $C_g=C_l\Pi^{+}$ (notice that $\Pi^{+}=T^{-1}\begin{bmatrix} I_l \\ 0 \end{bmatrix}$) and hence $\hat{C}_s=C_s\Pi^{+}$.
 Define now
 \( \mathscr{I}_{\infty}(\hat{v}_0,g) = \limsup_{t \rightarrow \infty} \mathscr{I}_t(\hat{v}_0,g) \).
  Recall from \eqref{lpg:eq:cost1} and \eqref{lpg:eq:cost2} the definition of the cost functions $\mathscr{J}_{\infty}$ and $\mathscr{J}_t$.
 It is not hard to see that:
 \begin{equation}
 \label{pf:main:eq1}
 \begin{split}
 & \mathscr{I}_t(\Pi(v_0),g)=\mathscr{J}_t(v_0,g)\,, \\
 & \mathscr{I}_{\infty}(\Pi(v_0),g)=\mathscr{J}_{\infty}(v_0,g),
 \end{split}
 \end{equation}
 for any initial state $v_0$ of $\mathscr{S}$ such that $v_0 \in \V_g$.

 Consider now the problem of minimizing $\lim_{t \rightarrow +\infty} \mathscr{I}_t(\hat{v}_0,g)$. Let us apply the feedback transformation $g=\hat{F}v+Uw$ to $\mathscr{S}_g=(A_g,B_g,C_g,D_g)$, where $\hat{F}$ and $U$ are defined as follows.
 If $D_l$ is of full column rank, then $U=(D^T_gSD_g)^{-1/2}$ and $\hat{F}=-(D^T_gSD_g)^{-1}D^T_gSC_g$.
  %This corresponds to the transformation
  %described in \cite[Section 10.5, eq. (10.32)]{TrentelmanBook}.
  If $D_l=0$, then $\hat{F}=0$ and $U=1$.
  Recall the ODE-LTI \eqref{opt:st2} and the corresponding cost function $\widehat{\mathscr{I}}_t(\hat{v}_0,w)$
  from the proof of Lemma \ref{opt:control:lemma1}.
   By construction of $\hat{F}$ and $U$,
  $\|S^{1/2}(C_g+D_g\hat{F})v(t)+D_gUw(t)\|^2=v^T(t)(C_g+D_g\hat{F})^TS(C_g+D_g\hat{F})v(t)+w^T(t)w(t)$.
 Using these remarks, it is then easy to see that  $\mathscr{I}_t(\hat{v}_0,g)=\widehat{\mathscr{I}}_t(\hat{v}_0,w)$ for $g=\hat{F}v+Uw$.

   Consider now the problem of minimizing $\lim_{t \rightarrow \infty} \widehat{\mathscr{I}}_t(\hat{v}_0,w)$.
   First, we assume that $D_g=D_l$ is full column rank.
   We apply \cite[Theorem 3.7]{KwakernaakBook}. In the proof of Lemma \ref{opt:control:lemma1} it was already shown that
   $(A_g+B_g\hat{F},B_gU)$ is stabilizable and $(S^{1/2}(C_g+D_g\hat{F}),A_g+B_g\hat{F})$ is observable.
   Let us now return to the minimization problem. Notice that \eqref{ARE} is equivalent to the algebraic Riccati equation
   described in \cite[Theorem 3.7]{KwakernaakBook} for the problem of minimizing
   $\lim_{t \rightarrow \infty} \widehat{\mathscr{I}}_t(\hat{v}_0,w)$.
 From Lemma \ref{opt:control:lemma1} it follows that \eqref{ARE}  has a unique positive definite solution $P$,
  and $A_g+B_g\hat{F}-B_gUU^TB_g^TP=A_g-B_gK$ is a stable matrix.
   From \cite[Theorem 3.7]{KwakernaakBook}, there exists $w^{*} \in L^2([0,+\infty),\mathbb{R}^k)$ such that
   $\lim_{t \rightarrow \infty} \widehat{\mathscr{I}}_t(\hat{v}_0,w^{*})$ is minimal and $\hat{v}_0^TP\hat{v}_0 = \lim_{t \rightarrow \infty} \widehat{\mathscr{I}}_t(\hat{v}_0,w^{*})$. On the other hand, \cite[Theorem 3.7]{KwakernaakBook} also implies that $\hat{v}^T_0P\hat{v}_0=\lim_{t \rightarrow \infty} \inf_{w \in L^1_{loc}([0,t],\mathbb{R}^k)} \widehat{\mathscr{I}}_{t}(\hat{v}_0,w)$. Hence, $g^{*}=\hat{F}v^{*}+Uw^{*}$ satisfies \[
 \mathscr{I}_{\infty}(\hat{v}_0,g^*)=\hat{v}_0^TP\hat{v}_0 = \lim_{t \rightarrow \infty}\inf_{g \in L^1_{loc}([0,t],\mathbb{R}^k)} \mathscr{I}_t(\hat{v}_0,g)\,,\] where $\dot v^* = (A_g+B_g\hat{F})v^*+B_gUw^*$, $v^*(0)=\hat{v}_0$.
  A routine computation reveals that $(v^{*},g^{*})$ satisfies
 \[ \dot v^* = A_gv^{*}+B_gg^* \mbox{ and } g^*=-Kv^* \mbox{ and } v^*(0)=\hat{v}_0. \]

Assume now that $D_l=D_g=0$. Then $B_g=0$ and $k=1$. Moreover, in this case, the solution $P=P^T > 0$ of \eqref{ARE} is in fact the observability grammian of the ODE-LTI
$\dot v = A_gv$, $y=S^{1/2}C_gv$. From the well-known properties of observability grammian it then follows that
$\hat{v}_0^TP\hat{v}_0=\int_0^{\infty} y^T(s)y(s)ds = \lim_{t \rightarrow \infty} \widehat{\mathscr{I}}_t(\hat{v}_0,0)=\widehat{\mathscr{I}}_{\infty}(\hat{v}_0,0)$.
Since $B_g=0$, $\widehat{\mathscr{I}}_t(\hat{v}_0,0)= \widehat{\mathscr{I}}_t(\hat{v}_0,w)=\mathscr{I}_t(\hat{v }_0,w)$ ($g=Fv+Uw=w$ if $D_g=0$) for any
 $w \in L^1_{loc}([0,+\infty),\mathbb{R})$, $t \in [0,+\infty)$, hence $\hat{v}_0^TP\hat{v}_0  =  \lim_{t \rightarrow \infty}\inf_{g \in L^1_{loc}([0,t],\mathbb{R}^k)} \mathscr{I}_t(\hat{v}_0,g)\).

Hence,  in both cases ($D_g=D_l$ is full rank or $D_l=0$),
from Theorem \ref{lqg:eq} it then follows that $({x^*}^T,{u^*}^T)^T=C_gv^*+D_gg^*$ is a solution of the infinite horizon optimal control problem and that $(x^*,u^*)$ satisfies \eqref{opt8} and \eqref{opt9}.

\textbf{(i) $\implies$ (iii)} If $(x^*,u^*) \in \mathscr{D}_{z}(\infty)$
  is a solution of the infinite horizon optimal control problem, then, by definition,
  $+\infty > J_{\infty}(x^*,u^*)=J^{*}=\limsup_{t \rightarrow \infty}  \inf_{(x,u) \in \mathscr{D}_{z}(t_1))} J_t(x,u)$.

\textbf{(iii) $\implies$ (ii)}
 From \eqref{daelincost} it follows that the condition of \textbf{(iii)}
 implies that there exists $M > 0$ such that for all $t > 0$,
 \begin{equation}
 \label{minimax:pf1}
  \inf_{g \in L^1([0,t],\mathbb{R}^k)} \mathscr{J}_t(\LMAP(z),g) \le M.
 \end{equation}
 Recall that
 $S=\mathrm{diag}(Q,R)$ and
 for each $g \in L^1_{loc}(I,\mathbb{R}^k)$, $[0,t] \subseteq I$, and initial state $v_0$, define
 \[
   \begin{split}
    \mathscr{H}_t(v_0,g) = \int_0^{t} (C_lv(s)+D_lg(s))^TS(C_lv(s)+D_lg(s))ds   \\
    \dot v = A_lv+B_lg \mbox{, \ } v(0)=v_0\,.
   \end{split}
 \]
 It then follows that $\mathscr{H}_t(v_0,g) \le \mathscr{J}_t(v_0,g)$ for any
 $t > 0$ and $\mathscr{H}_t(v_0,g)$ is non-decreasing in $t$.
 Hence, \eqref{minimax:pf1} implies that
 \[ \forall t \in (0,+\infty): \inf_{g \in L^1_{loc}([0,t],\mathbb{R}^k)} \mathscr{H}(v_0,g) < M. \]
 From classical linear theory \cite{TrentelmanBook,KwakernaakBook}
 it follows that if $H$ is the unique symmetric, positive semi-definite solution of the Riccati equation
 \begin{equation}
  \label{minimax:pf:DARE}
  \begin{split}
     & \dot H(t)= A_l^TH(t)+H(t)A_l-K^T(t)(D_l^TSD_l)K(t)+C_l^TSC_l  \\
     & K(t)=(D_l^TSD_l)^{-1}(B_l^TH(t)+D_l^TSC_l) \mbox{ and } H(0)=0 \\
\\
  \end{split}
  \end{equation}
 then $v^T_0H(t)v_0=\inf_{g \in L^1_{loc}([0,t],\mathbb{R}^k)} \mathscr{H}_t(v_0,g)$. The latter may be easily seen applying the state feedback transformation $g=\hat{F}v+Uw$ with $F,U$ defined as in the proof of Theorem~\ref{finhopr:opt} and solve the resulting standard LQ control problem for the transformed system. It then follows that \eqref{minimax:pf:DARE} is the differential
Riccati equation which is associated with this problem.
%  \footnote{One way to see this is to apply the state feedback transformation
%  $g=\hat{F}v+Uw$, where $U=-(D^T_lSD_l)^{-1/2}$ and $\hat{F}=-(D^T_lSD_l)^{-1}D^T_lSC_l$ and solve the resulting standard LQ control problem for the transformed system. Namely, in this case the cost transforms
% $\mathscr{I}(v_0,g,t_1)=\int_0^{t_1} (v^T(t)(C_g+D_g\hat{F})^TS(C_g+D_g\hat{F})v(t)+w^T(t)w(t))dt$ where $\dot v = A_lv+B_lg$.
% It then follows that \eqref{minimax:pf:DARE} is the differential
% Riccati equation which is associated with the problem of
% minimizing
%  $\int_0^{t_1} (v^T(t)(C_g+D_g\hat{F})^TS(C_g+D_g\hat{F})v(t)+w^T(t)w(t)dt$.}
 Note that the matrix $\dot H(t)$ is symmetric and positive semi-definite,
 since $b^TH(t)b$ is monotonically non-decreasing
 %{\bf Misha, monotonically  non-decreasing implies monotonically increasing, isn't it? it implies that b^TH(t)b \le b^TH(s)b for t \le s, i.e. it i monotonically non-strictly increasing. However, monotonically increasing is sometimes used in the sense of strictly monotonically increasing, hence the terminology}
for all $b$.
%Indeed, symmetry follows immediately from~\eqref{minimax:pf:DARE} and semi-definiteness may be proved as follows: for every state $v_0$,
 %$v_0^TH(t)v_0$ is a non-decreasing differentiable function of $t$ and hence its derivative $v_0^T\dot H(t)v_0$ is non-negative.
 Define the set
 \begin{equation}
 \label{minimax:pf:set}
     V=\{ v_0 \mid \sup_{t \in (0,+\infty)} v_0^TH(t)v_0 < +\infty \}.
 \end{equation}
  By assumption \textbf{(iii)} it follows that $\LMAP(z) \in V$. From \cite[Theorem 10.13]{TrentelmanBook} it follows that\footnote{Indeed, recall the system $\mathscr{S}_g=(A_g,B_g,C_g,D_g)$ and the map $\Pi$. Recall that $\mathscr{S}_g$ is stabilizable and hence by~\cite[Theorem 10.19]{TrentelmanBook} for every $v_0 \in \V_g$ there exists an input $g \in L^1([0,+\infty),\mathbb{R}^k)$ such that $\int_0^{\infty} ((C_gv(s)+D_gg(s))^TS(C_gv(s)+D_gg(s))ds < +\infty$ for $\dot v = A_gv+B_gg$, $v(0)=\Pi(v_0)$.
Consider the state trajectory $r(t)$ where $\dot r=A_lr+B_lg$, $r(0)=v_0$.
It then follows that $\Pi(r)=v$ and hence $C_lr+D_lg=C_gv+D_gg$.
 Therefore,
$v_0^TH(t_1)v_0 \le \mathscr{I}_t(v_0,g) \le \int_0^{\infty}((C_gv(s)+D_gg(s))^TS(C_gv(s)+D_gg(s))ds < +\infty$.}
$\V_g \subseteq V$. It is also easy to see that $V$ is a linear space.
  We will show that $V=\V_g$, and so $\LMAP(z) \in \V_g$ follows.

  First, we will argue that $V$ is invariant with respect to $A_l$.
  To this end, consider $v_0 \in V$ and set $v_1=e^{-A_lt}v_0$.
  For any $t_1 \in [0,+\infty)$ and any $g \in L^1_{loc}([0,t_1],\mathbb{R}^k)$, define
  $\hat{g} \in L^1_{loc}([0,t_1+t),\mathbb{R}^k)$ as $\hat{g}(s)=0$, $s \le t$ and
  $\hat{g}(s)=g(s-t)$ if $s > t$.
  Consider $\dot v = A_lv + B_l\hat{g}$, $v(0)=v_1$. It then follows
  $v(t)=v_0$ and hence
  %that $v(s)=e^{A_st}v_1$ for $s \le t$,
  %and if $\dot r = A_lr+B_lg$, $r(0)=v_0$,
  %then $v(t+s)=r(s)$ for all $s \in [0,+\infty)$.
  \begin{equation}
 \label{minimax:pf:set1}
 \begin{split}
  & \mathscr{H}_{t+t_1}(v_1,\hat{g})=  \\
      & \mathscr{H}_{t_1}(v_0,g) + \int_0^{t} (C_le^{A_ls}v_1)^TS(C_le^{A_ls}v_1)ds
   \end{split}
 \end{equation}
  Since $v_0 \in V$, it then follows that there exists $\Gamma > 0$ such that
   for any $t_1$ there exists $g \in L^{2}([0,t_1],\mathbb{R}^k)$ such that
  $v_0^TH(t_1)v_0 = \mathscr{H}_{t_1}(v_0,g) \le \Gamma$. Hence, from
  \eqref{minimax:pf:set1} it follows that
  \(
    %%Note that \hat{g} denotes the inputs which are zero before $t$. However, the infimum has to be taken over all the input \bar{g}.
      v_1^TH(t_1+t)v_1 = \inf_{\bar{g} \in L^1_{loc}([0,t_1+t],\mathbb{R}^k)} \mathscr{H}_{t+t_1}(v_1,\bar{g}) \le
     \Gamma+\int_0^{t} (C_le^{A_ls}v_1)^TS(C_le^{A_ls}v_1)ds < +\infty,
  \)
 and hence
 \[
     \sup_{t_1 \in (0,+\infty)} v_1^TH(t_1)v_1  \le
    \sup_{t_1 \in (0,+\infty)} v_1^TH(t_1+t)v_1 < +\infty.
  \]
In the last step we used that
  $v_1^TH(t+t_1)v_1 \ge v_1^TH(t_1)v_1$ for all $t,t_1 \in [0,+\infty)$.
 Hence, $v_1 = e^{-A_lt}v_0 \in V$. Since $V$ is a linear space,
and $t$ is arbitrary, it then follows that
$A_lv_0=-\frac{d}{dt} e^{-A_lt}v_0|_{t=0} \in V$. That is,
 $V$ is $A_l$ invariant.

Notice that the controllability subspace of $\mathscr{S}$ is contained
in $\V_g \subseteq V$ and that $\IM B_l$ is contained in the controllability subspace $\mathscr{S}$. Hence, $\IM B_l \subseteq V$.
Now, for any $b \in V$, the function $b^TH(t)b$ is monotonically non-decreasing in $t$ and it is bounded, hence $\lim_{t \rightarrow \infty} b^TH(t)b$ exists and it is finite. Notice for any $b_1,b_2 \in V$, $(b_1+b_2)^TH(t)(b_1+b_2)=b^T_1H(t)b_1+2b_1^TH(t)b_2+b_2^TH(t)b_2$ and as $b_1+b_2 \in V$, the limit on both sides exists and so the limit $\lim_{t \rightarrow \infty} b_1^TH(t)b_2$ exists. Consider now a basis $b_1,\ldots,b_r$ of $V$ and for any $t$ define $\hat{H}_{i,j}(t)=b_i^TH(t)b_j$. It then follows that the matrix $\hat{H}(t)=(H_{i,j}(t))_{i,j=1,\ldots,r}$ is
  positive semi-definite, symmetric and there exists a positive semi-definite
  matrix $\hat{H}_{+}$ such that $\hat{H}_{+}=\lim_{t \rightarrow \infty} \hat{H}(t)$.
  From \eqref{minimax:pf:DARE} and $A_lV \subseteq V$ it follows
  that $A_lb_i \in V$ for all $i=1,\ldots,r$ and hence
  $\lim_{t \rightarrow t} b_i^TH(t)A_lb_j$,
  $\lim_{t \rightarrow t} b_i^TA_l^TH(t)b_j$ exist for all $i,j=1,\ldots,r$.
  %Notice that the term
  %\( b^T_iK(t)^T(D_l^TSD_l)K(t)b_j \) is the sum of
  %the terms $b^T_i(B_l^TH(t))^T(D_l^TSD_l)^{-1}(B_l^TH(t))b_j$
  %  $(D_l^TSC_l b_i)^T(D_l^TSD_l)^{-1}(B_l^TH(t)b_j)$, $(D_l^TSC_l b_j)^T(D_l^TSD_l)^{-1}(B_l^TH(t)b_i)$ and the term \\
 %$(D_l^TSC_l b_i)^T(D_l^TSD_l)^{-1}(D_l^TSC_l b_j)$.
  From $\IM B_l \subseteq V$ and the fact that for any $x,z \in V$, the limit
  $\lim_{t \rightarrow \infty} x^TH(t)z$ exists, it follows that
  the limits $\lim_{t \rightarrow \infty} B_l^TH(t)x=\lim_{t
    \rightarrow \infty} x^TH(t)B_l$ exist for all $x \in V$. Applying this remark to $x=b_i$ and $x=b_j$,
 it follows that $\lim_{t \rightarrow \infty} b^T_iK(t)^T(D_l^TSD_l)K(t)b_j$ exists.
%{\bf Misha: why $(D_l^TSC_l b_i)^T(D_l^TSD_l)^{-1}(B_l^TH(t)) b_j$ exists? You did not prove that $x^T H(t) b_i$ converges for any $x$ which might be outside $V$.}
Hence, for any
  $i,j=1,\ldots,r$, the limit of $\dot{\hat{H}}_{i,j}(t)=b_i^T\dot H(t) b_j$ exists
  as $t \rightarrow \infty$ and hence the
  limit $\lim_{t \rightarrow \infty} \dot{\hat{H}}(t)=:Z$ exists. Moreover,
  since $\dot H(t)$ is symmetric and positive semi-definite, it follows that $Z$
  is symmetric and positive semi-definite.

  We claim that $Z$ is zero. To this end it is sufficient to show that $\lim_{t \rightarrow \infty} b^T\dot H(t)b = 0$ for any $b \in V$. Indeed, from this it follows that $Z_{i,j}=\lim_{t \rightarrow \infty} 0.5((b_i+b_j)^T\dot H(t)(b_i+b_j)-b^T_i\dot H(t)b_i-b_j^TH(t)b_j)=0$ for any $i,j=1,\ldots,r$. Now, fix $b \in V$ and assume that $c=\lim_{t \rightarrow \infty} b^T\dot H(t)b \ne 0$.
  Set $h(t)=b^TH(t)b$. It then follows that $\dot h(t)=b^T\dot H(t)b$ and thus there exists $T > 0$ such that for all
  $t > T$, $\dot h(t) > \frac{c}{2} > 0$ (recall that $\dot H(t)$ is positive semi-definite).
  Hence, $h(t)=h(0)+\int_0^{t} \dot h(s)ds > \int_{T}^{t} \dot h(s)ds > (t-T)\frac{c}{2}$.
  Hence, $h(t)$ is not bounded, which contradicts to the assumption that
  $b \in V$.

  Hence, $Z=0$ and thus, it follows that $\hat{H}_{+}$ satisfies
  the algebraic Riccati equation
  \begin{equation}
  \label{minimax:pf:ARE}
  \begin{split}
     & 0 = \hat{A}_l^T\hat{H}_{+}+\hat{H}_{+}\hat{A}_l-K^T(D_lSD_l)K+\hat{C}_l^TS\hat{C}_l   \\
     & K=(D_l^TSD_l)^{-1}(\hat{B}_l^T\hat{H}_{+}+D_l^TS\hat{C}_l) \\
  \end{split}
  \end{equation}
  where $\hat{A}_l,\hat{B}_l,\hat{C}_l$ are defined as follows:
  $\hat{A}_l$ and $\hat{C}_l$ are the matrix representations of
  the linear maps $A_l$ and $C_l$ restricted to $V$,
  $\hat{B}_l$ is the matrix representation of the map  $\mathbb{R}^k \ni g \mapsto B_lg \in V$ in the basis $b_1,\ldots,b_r$ of $V$ chosen as above.
  Note that for $\hat{A}_l$ and $\hat{B}_l$ to be well defined,
  we had to use the facts
  $\IM B_l \subseteq V$ and $A_l V \subseteq V$.
  Notice that $D_l$ is injective as a linear map and
  recall from Remark \ref{geo:rem1} that
  $C_lp+D_lg=0$ implies that $p=0$, $g=0$ and hence the largest output nulling
  subspace $\V(\hat{\Sigma})$
  of the linear system $\hat{\Sigma}=(\hat{A}_l,\hat{B}_l,\hat{C}_l,D_l)$ is zero.
  Then from \cite[Theorem 10.19]{TrentelmanBook}, $\V(\hat{\Sigma})=0$ and
  \eqref{minimax:pf:ARE} it follows that
  $\hat{\Sigma}$ is stabilizable. Since $\Sigma$ is just the restriction of
  $\mathscr{S}$ to $V$, it then follows that every state from $V$ is
  stabilizable and hence $V=\V_g$.

Let us now prove that the optimal solution satisfies \eqref{opt8},\eqref{opt9}, $u^*=K_fx^*$ and that it is a solution of \eqref{opt:feedback2} such that $Ex^{*}=z$.
   From the proof of \textbf{(i) $\rightarrow$ (ii)} it follows that $(x^*,u^*)$ from \eqref{opt8}
   is a solution of the infinite horizon optimal control problem for the initial state $z$ and that it satisfies \eqref{opt9}.
   Moreover, from \eqref{opt8} it follows that
   $u^*=K_fx^*$.
  Finally, all the solutions $(x,u)$ of \eqref{opt:feedback2}  on $[0,+\infty)$ satisfy
  $(x^T,u^T)^T(s)=(C_{g}-D_{g}K)v(s)$ a.e., $\dot v(s) = (A_g-B_gK)v(s)$, $v(s)=\LMAP_g(Ex(s))$, $s \in [0,+\infty)$.
   Hence, $(x^*,u^*)$ is indeed a solution of \eqref{opt:feedback2} such that  $Ex^*(0)=z$ and any other solution $(x,u)$ of \eqref{opt:feedback2} with $Ex(0)=z$ satisfies
   $x=x^*$, $u=u^*$ a.e.
\end{IEEEproof}

\section{Numerical example}
\label{sect:num}
 The spread of heat in a body can be described by the partial differential equation (PDE) of heat transfer,
 \begin{equation}
 \label{heat:eq1}
  \dfrac{dV(t)}{dt} = \mathcal{A}V(t) +Bu(t),
 \end{equation}
 where $V(t),u(t)$ take values in the space $H:=L^2([-1,1],\mathbb{R})$ and $\mathbb{R}^{N_u}$ respectively, and $\mathcal{A}v = -c^2\dfrac{d^2}{dx^2} v$ for all $v \in H$ such that $\mathcal{A}v\in H$, $Bv = \sum_{k=1}^{N_u} \sin(k\pi x) v_k$  for all $v=(v_1,\ldots,v_{N_u})^T \in \mathbb{R}^{N_u}$.
We say that a pair $(V,u)$ solves~\eqref{heat:eq1} if $u \in L^2_{loc}([0,+\infty), \mathbb{R}^{N_u})$, and
$V:[0,+\infty) \rightarrow H$ is a Frechet differentiable function such that $\mathcal{A}V(t)$ is defined, and \eqref{heat:eq1} holds a.e., and $V(t)(-1)=V(t)(1)=0$. The PDE \eqref{heat:eq1} is the simplest model for heat transfer. From a control perspective, $u$ is the control input. We would like to find $u$ with minimal energy that minimizes the energy of the heat distribution, i.e., \( \mathscr{J}(V,u)=\int_0^{\infty} \|V(t)\|_H^2 + \|u(t)\|^2 dt \) is minimal. Here $\|.\|_H$ is the standard norm of the space $H$, i.e.
$\|\phi\|_H^2=\int_{-1}^{1} \| \phi(\upsilon)\|^2_2 d\upsilon$. 

To compute the optimal control PDEs are usualy approximated by finite-dimensional models. Below we will do the same: we will present a DAE-LTI model which approximates \eqref{heat:eq1}, and, accordingly, the problem of minimizing $\mathscr{J}(V,u)$ is approximated by an instance of Problem \ref{opt:contr:def}, which cannot be solved by the existing methods. We stress that the equation \eqref{heat:eq1} has an eigen-basis $\{\sin(k\pi \upsilon)\}$ that allows to compute an exact finite-dimensional model representing the projection of~\eqref{heat:eq1} onto a given finite-dimensional subspace. This, in turn, allows to find an exact minimizer
$(V_{opt},u_{opt})$ of $\mathscr{J}$. We will compare this minimizer against the control law obtained by solving Problem \ref{opt:contr:def}. Note, that for more complicated PDEs, the eigen-basis may not be available and, hence, it may be impossible to compute the exact minimizer of $\mathscr{J}$. However, the method based on approximating a PDE by a DAE-LTI does not rely upon the existence of this basis and applies in the general case. Our experiment shows that the trajectory of the DAE-LTI approximates well the optimal solution of \eqref{heat:eq1}, and
the corresponding control law performs well when applied to \eqref{heat:eq1}. 
Note that this control law does not enforce the optimal solution of the DAE-LTI, when applied to the original model it enforces a near optimal solution.
The practical signficance of this is as follows. The DAE-LTI model and the corresponding instance of Problem \ref{opt:contr:def} can be applied to PDEs which are more general than
\eqref{heat:eq1}. However, the exact optimum of $\mathscr{J}(V,u)$ can be calculated only for \eqref{heat:eq1}, even small modifications of \eqref{heat:eq1} to render it more realistic
destroy this property. That is, DAE-LTIs could be a valuable tool for controlling systems described by PDEs. The reason we chose a model for which the control problem can be solved
exactly, is that it allows us to compare the perfomance of DAE-LTI based control with the true optimal one. Otherwise, it is not clear how to evaluate the performance of DAE-LTI based
controllers.

First we compute the optimal solution $V_{opt}$ and $u_{opt}$. To this end, consider  the orthogonal  basis $\{\sin(\pi k \upsilon)\}_{k=1}^{\infty}$ in $H$ and fix
an integer $N > 0$.  We will follow the classical projection Galerkin method, and represent the coordinates of the solution $V$ of \eqref{heat:eq1} as a solution of an ODE-LTI.
To this end, let $A_e$ be a diagonal matrix whose $i$th diagonal element is $-c^2i^2\pi^2$, $i=1,\ldots,N$, and let $B_e=I_{N_u}$.
Then for any $u \in L^2_{loc}([0,+\infty),\mathbb{R}^{N_u})$, $z(t)=(z_1(t),\ldots,z_N(t))^T$ is a solution of
\begin{equation}
\label{heat:eq1:exact}
   \dot z = A_{e}z+B_{e}z,
\end{equation}
if and only if
$(V,u)$, $V(t)(\upsilon)=\sum_{k=1}^{N} z_k(t) \sin(\pi k \upsilon)$ is a solution of \eqref{heat:eq1} and $\mathscr{J}(V,u)=J_e(z,u)\stackrel{def}{=} \int_0^{\infty} \|z(t)\|^2+\|u(t)\|^2 dt$.
Hence, if we want to minize $\mathscr{J}(V,f)$ with the additional restriction that $V(0)(\upsilon)=\lambda \sin(\pi n \upsilon)$ for some $n \le N$ and $\lambda > 0$, then it is equivalent to minimizing $J_{e}(z,u$ subject to \eqref{heat:eq1:exact} for $z(0)=\lambda e_n$. The latter problem is a standard LQ control problem. The optimal trajectory $(z_{opt},u_{opt})$ of \eqref{heat:eq1:exact} uniquely determines the optimal solution$V_{opt},u_{opt}$ of \eqref{heat:eq1}, and $u_{opt}=K_{opt}z_{opt}$ for some matrix $K_{opt}$ (equivalently, $u_{opt}=\mathcal{K}_{opt}V_{opt}$ for a suitable bounded operator $\mathcal{K}_{opt}$).  In particular, for $\lambda=10,n=34, N=40, N_u=35$ the optimal cost is $J_{e}(z_{opt},u_{opt})=\mathscr{J}(V_{opt},u_{opt})=3.94$. The approach presented above is a particular instance of the classical Galerkin projection method, applied to the basis $\{\sin(\pi n \upsilon)\}_{n=1}^{\infty}$. 

The approach above relies heavily on the fact that the functions $\{\sin(\pi n \upsilon)\}_{n=1}^{\infty}$ are eigenfunctions of $\mathcal{A}$, and hence cannot easily be generalized to other PDEs.However, it is difficult to find eigenfunctions for $\mathcal{A}$ it it is slightly modified so that it is more realistic. However, if we choose
just any basis and we apply the classical Galerkin projection method, then the resulting ODE-LTI will give a poor approximation of solutions of \eqref{heat:eq1}.
For this reason, we use the approach of \cite{ZhukCDC13} for constructing a DAE-LTIs whose solutions  approximate the solutions of \eqref{heat:eq1}. Consider the (non-orthogonal) basis $\{\phi_k\}_{k=0}^{\infty}$ of $H$, where $\phi_k=P_{k+1}-P_{k}$, $k \ge 1$, with $P_i$ denoting the $i$th Legendre polynomial, $i=0,1,\ldots$. In practice, this basis is often chosen due to its numerical properties. % when it is not clear which basis to choose.
Define the linear map $P_N$, $P_N^{+}$, and the matrices $M_N,\hat{M}_N, A_N$ as
\begin{equation*}
%\label{heat:dae:eq1}
  \begin{split}
  & \forall \psi \in H: P_N(\psi)=\begin{bmatrix} \langle\psi,\phi_1\rangle, & \ldots, & \langle\psi,\phi_N\rangle \end{bmatrix}^T \\
  & \forall z=(z_1,\ldots,z_N)^T \in \mathbb{R}^N: P_N^{+}(z)=\sum_{k=1}^{N} z_k \phi_k \\
   & \hat{M}_N=(\langle\phi_i,\phi_j\rangle)_{i,j=1}^{N}, ~ A_N=(\langle\mathcal{A} \phi_i,\phi_j\rangle)_{i,j=1}^{N} \\
   & M_N = \Lambda_N (\langle\phi_i,\phi_j\rangle)_{i,j=1}^{N}, ~ \Lambda_N=(\frac{2i+1}{2} \delta_{i,j})_{i,j=1}^{N} \\
  \end{split}
 \end{equation*}
 where $\langle\cdot,\cdot\rangle$ is the standard inner product in $H$ and $\delta_{i,j}$ is the Kronecker delta symbol for all $i,j=1,2,\ldots$. Intuitively, if 
 $\phi \in H$, then  $a=\hat{M}_N^{-1}P_N\phi$ is the vector of the first $N$ coordinates of $\phi$ w.r.t. $\{\phi_k\}_{k=1}^{\infty}$.
From \cite{ZhukCDC13} it  follows that if $(V,u)$ is a solution of \eqref{heat:eq1}, then with $a=\hat{M}_{N}^{-1}P_N V$, and $e_m=P_N\mathcal{A}\mathbf{e}$, $\mathbf{e}=V-P_N^{+}P_NV$, and $x=(a^T,e_m^T)^T$,  the pair $(x,u)$ is a solution of \eqref{eq:DAE}, where
 \[
  \begin{split}
   & E=\begin{bmatrix} M_N,  & 0 \end{bmatrix}, \quad  A=\begin{bmatrix} \Lambda_N A_N,  &  \Lambda_N \end{bmatrix}, \\
   & B=\Lambda_N \begin{bmatrix} P_N(\sin(\pi \upsilon)), & \ldots, & P_N(\sin(\pi N_u  \upsilon)) \end{bmatrix}^T, \\
   \end{split}
\]
Since according to \cite[Example 7.2, p. 121]{SpectralMethod}, $\langle\phi_i,\phi_j\rangle=4\frac{(2i+1)^2}{(2i-1)(2i+3)}\delta_{i,j}-\frac{2}{2i+3}\delta_{i+2,j}-\frac{2}{2i-1}\delta_{i-2,j}$, and $\langle\mathcal{A}\phi,\phi_j\rangle=-2c(2i+1)\delta_{i,j}$ for all $i,j=1,\ldots,N$, it follows that $E,A$ can easily be computed, while $B$ can be computed using numerical integration. The matrices $M_N,\Lambda_N,A,B,E,\hat{M}_N, \Lambda_N$ are presented in the supplementary material of this report for
$N=40$, $Nu=35,\mu=0.01$ and $c=\frac1{30}$. The corresponding matrices and the code for generating them can be downloaded from \verb'http://sites.google.com/site/mihalypetreczky/'.
Intuitively, the component $a(t)$  of $x(t)$ represents the first $N$ coordinates of $V(t)$ w.r.t $\{\phi_k\}_{k=1}^{\infty}$.
and $e_m$ is a linear function of the error $\mathbf{e}$ of approximating $V$ by its projection to the first $N$ basis vectors.
The smaller $\mathbf{e}$ is, the closer the behavior of \eqref{eq:DAE} is to that of \eqref{heat:eq1}. Conversely, if $(x,u)$ is a solution of \eqref{eq:DAE}, and $a$ is the vector of the first $N$ components of $x$, then we can define $V(t)=P_N^{+}a(t)$ and view $(V,u)$ as an approximation of a solution of \eqref{heat:eq1}. Note that, in general, $(V,u)$ obtained in this way will not be a true solution of \eqref{heat:eq1}.  
Consider Problem \ref{opt:contr:def}, with
\[ Q=\mathrm{diag}(\hat{M}_N, \mu I_N), \quad R=I_{N_u}, \quad Q_0=0 
\]
with $\mu > 0$ and  $(E,A,B)$ defined above. The choise of $Q$ is motivated by the following observation.
If $(x,u)$ is the solution of \eqref{eq:DAE} arising from a solution $(V,u)$ of \eqref{heat:eq1}, 
and we define $V_{approx}(t)=P_N^{+}a(t)$, then $\mathscr{J}(V_{approx},u)=J_{\infty}(x,u)-\mu \int_0^{\infty} \| e_m^2(t)\|_{2}^2dt$. 
Hence, if $\mathbf{e}=V-V_{approx}$ and $e_m$ are small,  i.e. $V_{approx}$ is a good approximation of $V$, then 
$J(x,u)$ is close to $\mathscr{J}(V_{approx},u)$, and hence to $\mathscr{J}(V,u)$, and it is reasonable to minimize $J_{\infty}(x,u)$ instead of $\mathscr{J}(V,u)$.
%More precisely,
% \(
%    J_{\infty}(x,u) =
%    \mathscr{J}(V,u) + \int_0^{\infty} \|\mathbf{e}(t)\|^2_H-\mu \|e_m(t)\|^2)dt,
% \)
% Hence, if $\mathbf{e}$ and $e_m$ are small, then $\mathscr{J}(V,u)$ is close to $J_{\infty}(x,u)$ and

\begin{figure}
\centering
\includegraphics[scale=0.3]{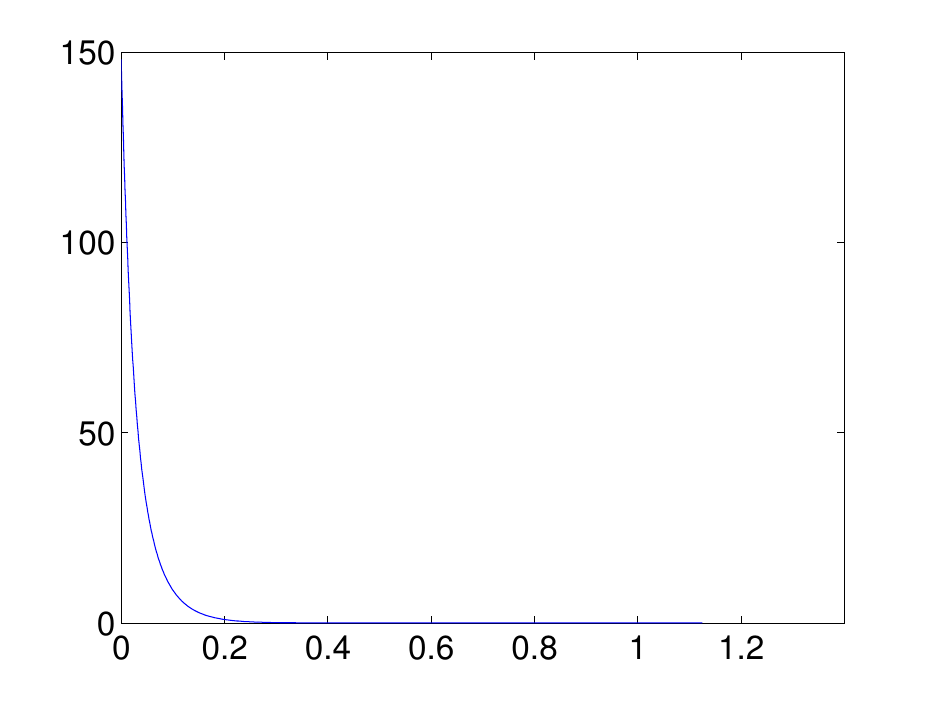}
\caption{Difference $e_{sol}(t)=\| V_{*}(t)- V_{opt}(t)\|_H^2 d\upsilon$ between the true optimum  and the optimal solution of \eqref{eq:DAE}\label{fig2}}
\end{figure}

\begin{figure}
\includegraphics[scale=0.4]{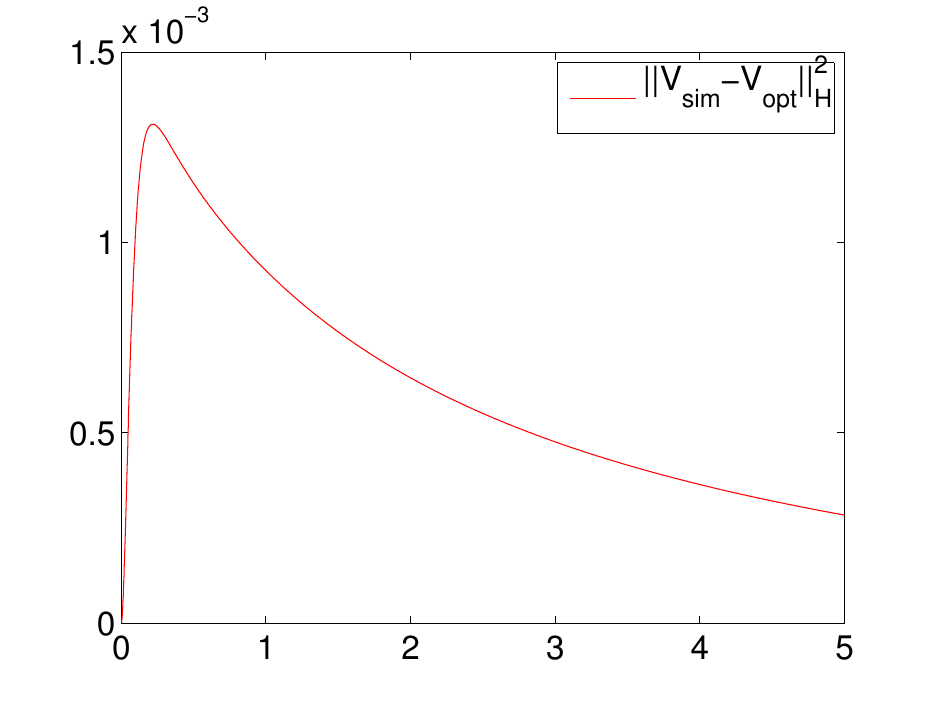}
\caption{Difference $e_{sim}(t)=\| V_{sim}(t)- V_{opt}(t)|_H^2 d\upsilon$ between the true optimal solution and the solution obtained by applying the feedback computed from the DAE-LTI \label{fig1}}
\end{figure}

The results of \cite{Kurina1,Kurina2,Kurina3,Kurina4,Stefanovski1,Stefanovski1.1, Stefanovski2,Stefanovski4} do not apply to this instance of Problem \ref{opt:contr:def}.
If it was the case, then there
would exist a feedback $u=Kx=K_1a+K_2e_m$, with the decomposition $x=(a^T,e_m^T)^T$,  %for \eqref{eq:DAE},
such that all the trajectories of the closed-loop system are at least stable. However, the closed-loop system reads as follows:
$\dot M_N a = (\Lambda_NA_N+BK_1)a+(\Lambda_N+BK_2)e_m$, and $e_m$ can be any square-integrable function. Hence, we can always choose $e_m$ so that it does not converge to zero as 
$t \rightarrow +\infty$. 
%$J_{\infty}(x,u)=+\infty$. This contradicts to the assumption that all the solutions of the closed-loop system are optimal.
This is not surprising, since $e_m$ is related to the error of approximating \eqref{heat:eq1} by a finite-dimnesional system, and this error cannot be controlled. \\
We computed a stabilizable associated ODE-LTI, solved \eqref{ARE}  and computed the optimal trajectory for the case:
$N=40$, $Nu=35,\mu=0.01$ and $c=\frac1{30}$. The corresponding matrices and the code for generating them and for solving the optimal control problem can be downloaded from \verb'http://sites.google.com/site/mihalypetreczky/', and they are also presented in the supplementary material of this report. From the specific initial condition $Ex(0)=\Lambda_N P_N V_{opt}(0)$,  % $V(0)(\upsilon)=\lambda \sin(\pi n \upsilon)$, $\lambda=10, n=34$,
 we computed the optimal trajectory $(x_*,u_*)$, and the optimal cost $J_{\infty}(x_*,u_*)=3.89$. As mentioned earlier, the solution $(V_{opt},u_{opt})$  of \eqref{heat:eq1} which satisfies $V_{opt}(0)(\upsilon)=\lambda \sin(\pi n \upsilon)$, $\lambda=10, n=34$
can be computed by using the eigen-basis of $\mathcal{A}$, see \cite{AutomaticaPaperArxive} for more
details. The true optimal value is $\mathscr{J}(V_{opt},u_{opt})=3.94$, which is close to $J_{\infty}(x_*,u_*)$. On Figure \ref{fig2}
we plotted the norm of the difference between $V_{*}(t)=(P_N^{+} a_*(t))$ and $V_{opt}(t)$, where $a_*(t)$ is the vector formed by the first $N$ components of $x_*$.
% $e_{sol}(t)=\int_{-1}^{1} \| V_{*}(t,\upsilon)- V_{opt}(t,\upsilon)|_2^2 d\upsilon$.
\\
Recall that the optimal solution $(x^*,u^*)$ satisfies the feedback $u^*=K_fEx^*$ for some matrix $K_f$.
Define now $\mathcal{K}_f=K_fM_N\hat{M}_N^{-1}P_N=K_f\Lambda_ NP_N$, and consider the feedback law $u= \mathcal{K}_f V$ for \eqref{heat:eq1}.
For a fixed initial value $V(0)$ the feedback $u=\mathcal{K}_fV$ yields a unique solution $(V,u)$ of \eqref{heat:eq1}, and $(V,u)$ yields a solution $(x,u)$ of
\eqref{eq:DAE} which satisfies $u=K_fEx$.
That is, some solutions of the closed-loop DAE-LTI \eqref{eq:DAE} with $u=K_fEx$ arise from the solutions of the closed-loop system
\eqref{heat:eq1} with the control law $u=\mathcal{K}_fV$. While the closed-loop DAE-LTI will have many solutions which minimize, the closed-loop PDE \eqref{heat:eq1}
with $u=\mathcal{K}_fV$ will have at most one.
We computed  a solution $(V_{sim},u_{sim})$, $u_{sim}=\mathcal{K}_fV_{sim}$ of \eqref{heat:eq1} on the interval $[0,T],T=5$, 
such that $V_{sim}(0)=V_{opt}(0)$. 
The cost  $\mathscr{J}(V_{sim},u_{sim})$ was approximated by $\mathscr{J}_T(V_{sim},u_{sim})=\int_0^T \| V_{sim}(t)\|_H^2 + \|u_{sim}(t)\|_{2}^2 dt$, yielding $3.96$ which is close to $\mathscr{J}(V_{opt},u_{opt})$.  Note that by increasing $T$, the value of
$\mathscr{J}_T(V_{sim},u_{sim})$ did not change significantly.  Furthermore, 
the difference between $V_{sim}(t)-V_{opt}(t)$  approaches zero as $t \rightarrow T$ and  and $max_{t \in [0,T]} \|V_{sim}(t,\upsilon)-V_{opt}(t,\upsilon)\|^2_Hd\upsilon = 0.0015$.  
The corresponding plot is presented in Figure \ref{fig1}.
Note that $V_{*} \ne V_{sim}$, since $(V_{*},u_*)$ need not be a true solution of \eqref{heat:eq1} while $(V_{sim},u_{sim})$ is so.
In fact, even the inital values are different: $V_{*}(0) \ne V_{opt}(0)=V_{sim}(0)$. This is not surprising, since $V_{*}(0)=P_N^{+}a_{*}(0)$, and 
$Ex_{*}(0)=\Lambda_N P_N V(0)$ implies that $a_*(0)$ is the vector of the first $N$ coordinates of $V_{opt}(0)$. However, since $V_{opt}(0)$ cannot be expressed as a 
finit sum of $\{\phi_k\}_{k=1}^{\infty}$, necessarily $P_{N}^{+} a_*(0) \ne V_{opt}(0)$. 

The solution $(V_{sim},u_{sim})$ was computed as follows. If $V_{opt}(0)(\upsilon)=\lambda \sin(n \pi \upsilon)$, then any solution $(V,u)$ of \eqref{heat:eq1} such that $V(0)=V_{opt}(0)$
satisfies $V(t)(\upsilon)=\sum_{k=1}^{N} z_k(t) \sin(k \pi \upsilon)$, where $(z,u)$ satisfies \eqref{heat:eq1:exact}. Hence, the feedback law $u=\mathcal{K}_fV=K_f E\hat{M}_N^{-1}P_N V$ is equivalent
to a feedback law $u=K_1 z$ for \eqref{heat:eq1:exact}, where the $j$th column of $K_1$ equals $K_fE\hat{M}_N^{-1}P_N \sin(\pi j \upsilon)$, $j=1,\ldots, N$.
We simulated
\eqref{heat:eq1:exact} on $[0,T]$, with $u=K_1z$ and the initial state being $z(0)=\lambda e_n$, where $e_n$ is the $n$th standard basis vector. 
For the example at hand, recall that $n=34$ and  $\lambda=10$, $T=5$.  
Let $z$ be the resulting state trajectory of \eqref{heat:eq1:exact}. We took $V_{sim}(t)=\sum_{k=1}^{N} z_k(t) \sin(k \pi \upsilon)$,  $u_{sim}(t)=K_1z(t)$ for all $t  \in [0,T]$.
It then follows that $(V_{sim},u_{sim})$ is the unique solution of \eqref{heat:eq1} such that $V_{sim}(0)=V_{opt}(0)$, and it satsfies $u_{sim}= \mathcal{K}_f V$. 
Moreover, $\mathscr{J}_T(V_{sim},u_{sim})=\int_0^T \|z(t)\|^2_2+\|u(t)|_2^2dt$.
The Matlab code for computing $V_{opt}$, $u_{opt}$, $V_*$,  $(x_*,u_*)$, $V_{sim},u_{sim}$, the corresponding cost functions and generating the plots can be downloaded from \verb'http://sites.google.com/site/mihalypetreczky/'. 

\begin{figure}
\includegraphics[scale=0.4]{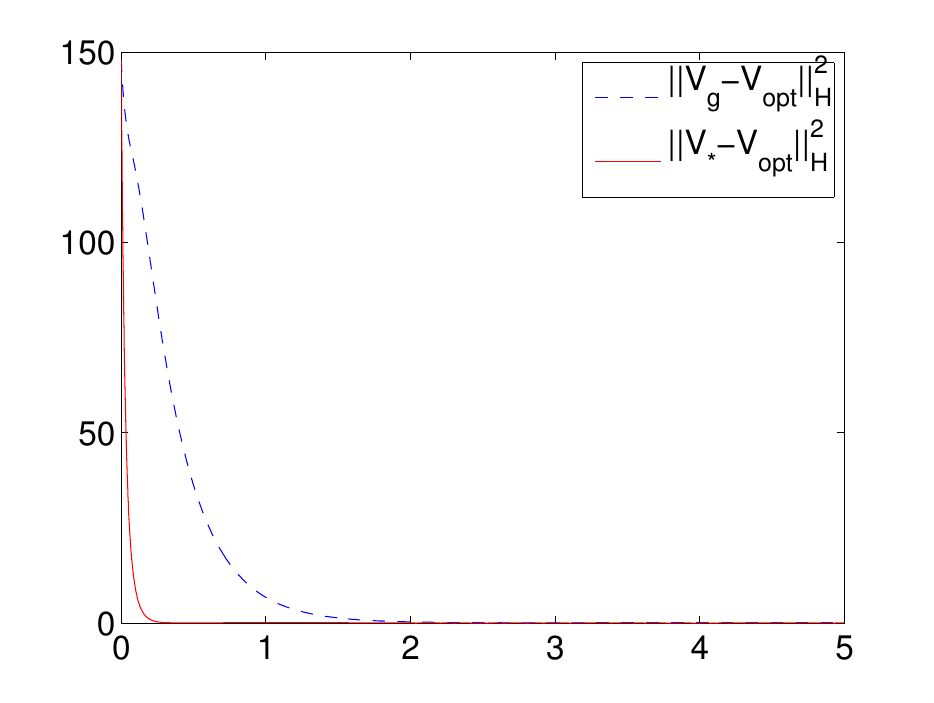}
\caption{Comparison with DAE-LTI: the difference $e_{g}(t)=\| (V_g(t)- V_{opt}(t)\|_H^2 d\upsilon$ between the true optimum  and the optimal solution of the ODE-LTI obtained via Galerkin method \label{fig3}, and $e_{sol}(t)=\| V_{*}(t)- V_{opt}(t)\|_H^2 d\upsilon$ }
\end{figure}

\begin{figure}
\includegraphics[scale=0.4]{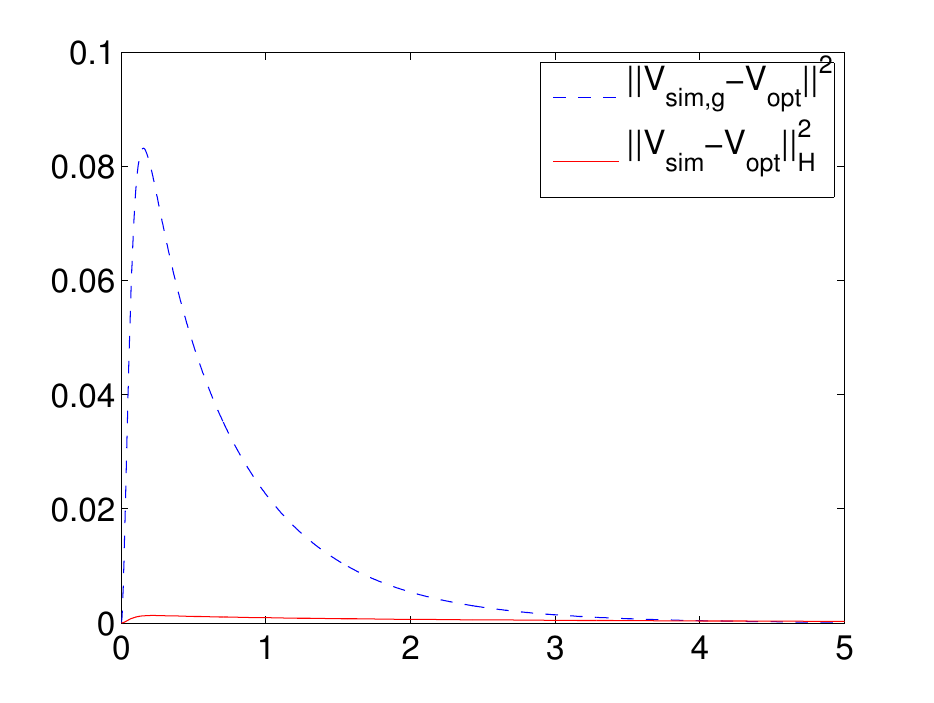}
\caption{Comparison with DAE-LTI: the difference $e_{sim,g}(t)=\| V_{sim,g}(t)- V_{opt}(t)|_H^2 d\upsilon$ between the true optimal solution and the solution obtained by applying the feedback computed from the ODE-LTI \label{fig4} and
$e_{sim}(t)=\| V_{sim}(t)- V_{opt}(t)|_H^2 d\upsilon$. 
 }
\end{figure}

Note that we could have tried to approximate solutions of \eqref{heat:eq1} using classical Galerkin projection method w.r.t. $\{\phi_{k}\}_{k=1}^{\infty}$. Then the problem of minimizing $\mathscr{J}(V,u)$ is approximated by a classical LQ control problem for an ODE-LTI. However,the estimate of the optimal value of $\mathscr{J}(V,u)$ and the perfomance of the corresponding controller obtained for the ODE-LTI are much worse than the one described in the previous paragraphs.
More precisely, in this case we approximate the first $N$ coordinates the solution $V$ of \eqref{heat:eq1} w.r.t. $\{\phi_k\}_{k=1}^{\infty}$ by the solution of
of the ODE-LTI $\dot a = M_N^{-1}A_Na+M_N^{-1} Bu$, and we approximate $\mathscr{J}(V,u)$ by $J_{g}(a,u)=\int_0^{\infty} a^T(t) \hat{M}_N a(t) dt + u^T(t)u(t)dt$. 
Instead of minimizing $\mathscr{J}(V,u)$ we then minimize $J_g(a,u)$. The latter is a classical LQ problem.
If $(a_g,u_g)$ is the minimal value of $J_g(a,u)$ subject to $\dot a = M_N^{-1}A_Na+M_N^{-1} Bu$, then $(V_g,u_g)$, with
$V_g(t)=P_N^{+} a_g$ can be viewed as an approximation of a solution of \eqref{heat:eq1}, and $\mathscr{J}(V_g,u_g)=J_g(a_g,u_g)$.
%%as this ODE-LTI does not take into account the approximation error.
For $a(0)=P_N V_{opt}(0)$, we the obtained estimate of $\mathscr{J}(V_{g},u_{g})=J_g(V_g,u_g)=6.13$, which is much higher than the true
value $\mathscr{J}(V_{opt},u_{opt})=3.94$. The plot of $\|V_g(t)-V_{opt}(t)\|_H^2$ and together with $\|V_{*}(t)-V_{opt}(t)\|_H^2$ is shown on Figure \ref{fig3}. It is clear that $V_g$ is much further from $V_{opt}$ than $V_{*}$ and it converges to it
much slower. 
When optimizing $J_g(q,u)$ we also obtain a feedback law $u=K_ga$. If can interpret the feedback $u=K_ga$ as a feedback $u=\mathcal{K}_g V$, $\mathcal{K}_g=K_g\hat{M}_N^{-1}P_N V$, and
apply it to \eqref{heat:eq1}. The resulting solution  $(V_{sim,g},u_{sim,g})$ of the closed-loop system is far from
the optimal $V_{opt}$ for $V_{sim,g}(0)=V_{opt}(0)$, and it converges slower to $V_{opt}$ than $V_{sim}$ obtained by applying the feedback $u=\mathcal{K}_fV$ computed by using DAE-LTI, see Figure \ref{fig4}. In fact, $\mathscr{J}_T(V_{simg,g},u_{sim,g})=5.55$, where $\mathscr{J}_T(V_{sim,g},u_{sim,g})=\int_0^T \| V_{sim,g}(t)\|_H^2 + \|u_{sim,g}(t)\|_{2}^2 dt \le \mathscr{J}(V_{sim,g},u_{sim,g})$. 
Hence, the value of $\mathscr{J}(V,{sim,g},u_{sim,g})$, is much larger than 
the optimum $\mathscr{J}(V_{opt},u_{opt})$. 
This is not surprising,
the approximation error, hence the resulting model is much less accurate. This example thus illustrates the advantages of using
a DAE-LTI model over the traditional Galerkin projection method.

In order to compute $(V_{sim,g},u_{sim,g})$ we used \eqref{heat:eq1:exact} in a similar manner as it was done to compute $(V_{sim},u_{sim})$.  That is, 
we used the fact that 
$V_{opt}(0)(\upsilon)=\lambda \sin(n \pi \upsilon)$, then any solution $(V,u)$ of \eqref{heat:eq1} such that $V(0)=V_{opt}(0)$
satisfies $V(t)(\upsilon)=\sum_{k=1}^{N} z_k(t) \sin(k \pi \upsilon)$, where $(z,u)$ satisfies \eqref{heat:eq1:exact}. Hence, the feedback law $u=\mathcal{K}_gV=K_g\hat{M}_N^{-1}P_N V$ is equivalent
to a feedback law $u=K_{1,g} z$ for \eqref{heat:eq1:exact}, where the $j$th column of $K_{1,g}$ equals $K_g\hat{M}_N^{-1}P_N \sin(\pi j \upsilon)$, $j=1,\ldots, N$.
We simulated
\eqref{heat:eq1:exact} on $[0,T]$, with $u=K_{1,g}z$ and the initial state being $z(0)=\lambda e_n$, where $e_n$ is the $n$th standard basis vector. 
For the example at hand, recall that $n=34$ and  $\lambda=10$, $T=5$.  
Let $z$ be the resulting state trajectory of \eqref{heat:eq1:exact}. We took $V_{sim,g}(t)=\sum_{k=1}^{N} z_k(t) \sin(k \pi \upsilon)$,  $u_{sim}(t)=K_{1,g}z(t)$ for all $t  \in [0,T]$.
It then follows that $(V_{sim,g},u_{sim,g})$ is the unique solution of \eqref{heat:eq1}, such that $V_{sim,g}(0)=V_{opt}(0)$, and it satsfies $u_{sim,g}= \mathcal{K}_f V$. 
Moreover, $\mathscr{J}_T(V_{sim,g},u_{sim,g})=\int_0^T \|z(t)\|^2_2+\|u(t)|_2^2dt$.

The code for computing $a_{g},u_g$, $V_g$, $V_{sim,g},u_{sim,g}$, the corresponding values of the cost functions and generating the plots can be downloaded from \verb'http://sites.google.com/site/mihalypetreczky/'. 

\section{Conclusions}
\label{sect:concl}
 We have presented a framework for representing solutions of non-regular DAE-LTIs as outputs of ODE-LTIs and for solving an LQ infinite horizon
 optimal control problem for non-regular DAE-LTIs. The solution concept we adopted for DAE-LTIs  allows non-differentiable solutions but excludes
 distributions.  We have shown that ODE-LTIs representing solutions of the same DAE-LTI are feedback equivalent and minimal. Moreover, we presented
 necessary and sufficient conditions for existence of a solution of the LQ infinite horizon control problem and an algorithm for computing the
 solution. The solution can be generated by  adding an algebraic constraint to the original DAE-LTI  and hence it can be interpreted as a result applying
 a controller in the behavioral sense \cite{WillemsInterconnect}.

 %We have presented the solution to infinite horizon
 %linear quadratic control problem for generic linear
 %DAEs. Specifically, for any tuple $(E,A,B)$ and initial condition
 %$x_0$ we proved that the
 %infinite horizon control problem is solvable if and only if the DAE
 %$\Sigma$ is stabilizable from $x_0$. We also proposed an intuitive
 %and computationally feasible way to check the stabilizability and
 %construct the optimal control which makes the proposed framework attractive for applications.

\bibliographystyle{plain}
%\bibliography{dae_bib}

%\appendix
\section{Matrices of the numerical example}
\label{sect:num:app}
 Please see the supplementary material of the report for the
 matrices of the numerical example. 
%\includepdf[fitpaper]{model_matrices1.pdf}
%\includepdf[fitpaper]{model_matrices2.pdf}
%\includepdf[fitpaper]{model_matrices3.pdf}
%\includepdf[fitpaper]{model_matrices4.pdf}
%\includepdf[fitpaper]{model_matrices5.pdf}

\end{document}